\documentclass[12pt]{article}
\usepackage{amssymb}
\usepackage{amsxtra}
\usepackage{amstext}
\usepackage{latexsym}
\usepackage{amsmath}
\usepackage{dsfont} 
\usepackage{yfonts}
\usepackage{mathrsfs}

\DeclareMathAlphabet{\mathpzc}{OT1}{pzc}{m}{it}

\newtheorem{definition}{Definition}[section]
\newtheorem{theorem}[definition]{Theorem}
\newtheorem{lemma}[definition]{Lemma}
\newtheorem{corollary}[definition]{Corollary}
\newtheorem{note}[definition]{Note}
\newtheorem{assumption}[definition]{Assumption}
\newtheorem{notation}[definition]{Notation}

\def\K{\mathbb K}
\def\Z{\mathbb Z}

\newcommand{\SL}{\mathfrak{sl}_2}
\newcommand{\uq}{U_q(\mathfrak{sl}_2)}

\begin{document}

\title{\bf Bidiagonal Triples}
\author{Darren Funk-Neubauer}
\date{}

\maketitle

\begin{abstract}
\noindent
We introduce a linear algebraic object called a bidiagonal triple.  A bidiagonal triple consists of three diagonalizable linear transformations on a finite-dimensional vector space, each of which acts in a bidiagonal fashion on the eigenspaces of the other two.  The concept of bidiagonal triple is a generalization of the previously studied and similarly defined concept of bidiagonal pair.  We show that every bidiagonal pair extends to a bidiagonal triple, and we describe the sense in which this extension is unique.  In addition we generalize a number of theorems about bidiagonal pairs to the case of bidiagonal triples.  In particular we use the concept of a parameter array to classify bidiagonal triples up to isomorphism.  We also describe the close relationship between bidiagonal triples and the representation theory of the Lie algebra $\SL$ and the quantum algebra $\uq$.    \\


\end{abstract}

\section{Introduction}

\noindent
In this paper we introduce a linear algebraic object called a bidiagonal triple.  Roughly speaking, a bidiagonal triple is a triple of diagonalizable linear transformations on a finite-dimensional vector space, each of which acts in a bidiagonal fashion on the eigenspaces of the other two (see Definition \ref{def:bdtriple} for the precise definition).  This concept arose as a generalization of the closely related concept of bidiagonal pair (see Definition \ref{def:bdpair}).  The theory of bidiagonal pairs was developed in \cite{Funk-Neubauer13}.  Any two of the three transformations in a bidiagonal triple form a bidiagonal pair.  Thus, a bidiagonal triple can be thought of as three bidiagonal pairs interwoven together.  \\

The main purpose of this paper is to generalize the results about bidiagonal pairs in \cite{Funk-Neubauer13} to the case of bidiagonal triples.  First we show that every bidiagonal pair can be extended to a bidiagonal triple (see Theorem \ref{thm:extend}).  That is, we present a construction showing how to build a bidiagonal triple starting from a bidiagonal pair.  In Theorem \ref{thm:unique} we describe the sense in which this construction is unique.  Using this close connection between bidiagonal pairs and bidiagonal triples we then do the following.  We associate to each bidiagonal triple a sequence of scalars called a parameter array (see Definition \ref{def:parray}).  We use this concept of parameter array to classify bidiagonal triples up to isomorphism (see Theorem \ref{thm:class}).  The statement of this classification does not make it clear how to {\it construct} the bidiagonal triples in each isomorphism class.  Hence, we show how to construct bidiagonal triples using finite-dimensional modules for the Lie algebra $\SL$ and the quantum algebra $\uq$ (see Theorem \ref{thm:uqsl2bidiag}).  The finite-dimensional modules for $\SL$ and $\uq$ are well known (see Lemmas \ref{thm:sl2mods} and \ref{thm:uq2mods}).  \\  

Bidiagonal triples and pairs originally arose in the study of the well-known quantum algebras $U_q(\widehat{\mathfrak{sl}}_2)$ and $\uq$.  The discovery of the so called equitable presentations of these algebras was the initial motivation for defining a bidiagonal triple.  These equitable presentations were discovered during the attempt to classify a linear algebraic object called a tridiagonal pair.  See below for more information on the equitable presentations and tridiagonal pairs.  Thus, the importance of bidiagonal triples lies in the fact that they provide insight into the relationships between several closely connected algebraic objects.  Although this paper is the first to explicitly define bidiagonal triples, they appear implicitly in \cite{Benkart04, Funk-Neubauer07, Funk-Neubauer09, ItoTer073, ItoTerWang06, Ter13}.  In \cite{Benkart04} bidiagonal triples were involved in constructing irreducible $U_q(\widehat{\mathfrak{sl}}_2)$-modules from the Borel subalgebra of $U_q(\widehat{\mathfrak{sl}}_2)$.  In \cite{Funk-Neubauer07} bidiagonal triples were used in constructing $U_q(\widehat{\mathfrak{sl}}_2)$-modules from certain raising and lowering maps satisfying the $q$-Serre relations.  In \cite{Funk-Neubauer09} bidiagonal triples were present in using a certain type of tridiagonal pair to construct irreducible modules for the $q$-tetrahedron algebra.  See below for more information on the $q$-tetrahedron algebra.  In \cite{ItoTer073} bidiagonal triples were also present in using a certain type of tridiagonal pair to construct irreducible $U_q(\widehat{\mathfrak{sl}}_2)$-modules.  Bidiagonal triples were used in \cite{ItoTerWang06} to show that the generators from the equitable presentation of $\uq$ have invertible actions on each finite-dimensional $\uq$-module.  See \cite[Theorem 18.3]{Ter13} for another example of how bidiagonal triples appear in the representation theory of $\uq$. \\

We now discuss the equitable basis for $\SL$ and the equitable presentation of $\uq$ because they will be used in our construction of bidiagonal triples.  The word equitable is used here because each exhibit a $\Z_3$-cyclic symmetry not appearing in the standard basis for $\SL$ or the standard presentation of $\uq$.  \\

The equitable basis for $\SL$ was first introduced in \cite{HarTer07} as part of the study of the tetrahedron Lie algebra \cite[Definition 1.1]{HarTer07}.  The tetrahedron  algebra has been used in the ongoing investigation of tridiagonal pairs, and is closely related to a number of well-known Lie algebras including the Onsager algebra \cite[Proposition 4.7]{HarTer07}, the $\SL$ loop algebra \cite[Proposition 5.7]{HarTer07}, and the three point $\SL$ loop algebra \cite[Proposition 6.5]{HarTer07}.  For more information on the tetrahedron algebra see \cite{Benkart07, Elduque07, Hartwig07, ItoTerinpress4}.  For more information on the equitable basis for $\SL$ see \cite{Al-Najjar10, Benkart10}.  \\

The equitable presentation of $\uq$ was first introduced in \cite{ItoTerWang06}.  This presentation was used in the study of the $q$-tetrahedron algebra \cite{ItoTer072}.  The $q$-tetrahedron algebra has been used in the ongoing investigation of tridiagonal pairs, and is closely related to a number of well-known algebras including the $U_q(\mathfrak{sl}_2)$ loop algebra \cite[Proposition 8.3]{ItoTer072} and the positive part of $U_q(\widehat{\mathfrak{sl}}_2)$ \cite[Proposition 9.4]{ItoTer072}.  For more information on the $q$-tetrahedron algebra see \cite{Funk-Neubauer09, ItoTer071, ItoTerinpress1, ItoTerinpress2, Miki10}.  We note that there exists an equitable presentation for the quantum algebra $U_q(\mathfrak{g})$, where $\mathfrak{g}$ is any symmetrizable Kac-Moody algebra \cite{Terinpress}. \\   

We now offer some additional information on tridiagonal pairs because of their close connection to bidiagonal triples.  Tridiagonal pairs originally arose in algebraic combinatorics through the study of a combinatorial object called a P- and Q-polynomial association scheme \cite{ItoTanTer01}.  The precise definition of a tridiagonal pair is given in \cite[Definition 1.1]{ItoTanTer01}.  Tridiagonal pairs appear in a wide variety of mathematical subjects including representation theory, orthogonal polynomials, special functions, partially ordered sets, statistical mechanics, and classical mechanics.  See \cite[Introduction]{Funk-Neubauer13} for the appropriate references.  A major classification result involving tridiagonal pairs appears in \cite{ItoTer09}.  A certain special case of a tridiagonal pair, called a Leonard pair, has also been the subject of much research.  The Leonard pairs are classified up to isomorphism \cite[Theorem 1.9]{Ter011}.  See \cite{Ter042} for a survey of Leonard pairs.  \\

This paper is organized as follows.  In Section 2 we define a bidiagonal triple and develop a number of related definitions to be used throughout the paper.  In Section 3 we recall the equitable basis for $\SL$, the equitable presentation of $\uq$, and the finite-dimensional modules for these algebras.  Section 4 contains statements of the main results of the paper.  In Section 5 we develop the tools needed to prove our main results including some additional properties of bidiagonal triples.  Sections 6 through 8 contain the proofs of the main results stated in Section 4.   

\section{Bidiagonal Triples}

In this section we recall the definition of a bidiagonal pair, present the definition of a bidiagonal triple, and make some observations about these definitions. 

\noindent
\begin{notation}
\rm
Throughout this paper we adopt the following notation.  Let $d$ denote a nonnegative integer.  Let $\K$ denote a field.  Let $V$ denote a vector space over $\K$ with finite positive dimension.  For linear transformations $X: V \to V$ and $Y: V \to V$ we define $[X, Y] := XY - YX$. 
\end{notation}

We present the following lemma in order to motivate the definitions of bidiagonal pair and bidiagonal triple.

\begin{lemma}
\label{thm:raise}
Let $X:V\to V$ and $Y:V \to V$ denote linear transformations.  
\begin{enumerate}
\item[\rm (i)] Suppose that there exists an ordering $\lbrace Y_i \rbrace _{i=0}^d$ of the eigenspaces of $Y$ with \\ $X Y_i \subseteq Y_i + Y_{i+1} \, (0 \leq i \leq d)$, where $Y_{d+1} = 0$.  Then for $0 \leq i \leq d$, the restriction $[X, Y] |_{Y_i}$  maps $Y_i$ into $Y_{i+1}$. 
\item[\rm (ii)] Suppose that there exists an ordering $\lbrace Y_i \rbrace _{i=0}^d$ of the eigenspaces of $Y$ with \\ $X Y_i \subseteq Y_{i-1} + Y_i \, (0 \leq i \leq d)$, where $Y_{-1} = 0$.  Then for $0 \leq i \leq d$, the restriction $[X, Y] |_{Y_i}$  maps $Y_i$ into $Y_{i-1}$. 
\end{enumerate}
\end{lemma}

\noindent
{\it Proof:}  (i) See the proof of \cite[Lemma 2.1]{Funk-Neubauer13}. \\
(ii) Similar to (i).  
\hfill $\Box$ \\

\begin{definition}
\rm
\cite[Definition 2.2]{Funk-Neubauer13}
\label{def:bdpair}
A {\it bidiagonal pair (BD pair) on $V$} is an ordered pair of linear transformations $A:V \to V$ and $A' :V \to V$ that satisfy the following three conditions.
\begin{enumerate}
\item[\rm (i)] Each of $A,\, A'$ is diagonalizable.
\item[\rm (ii)] There exists an ordering $\lbrace V_i \rbrace _{i=0}^d$ (resp.~$\lbrace V'_i \rbrace _{i=0}^D$) of the eigenspaces of $A$ (resp.~$A'$) with 
\begin{align}
\label{eq:bp1}
A' V_i \subseteq V_i + V_{i+1} \qquad (0 \leq i \leq d), \\
\label{eq:bp2}
A V'_i \subseteq V'_i + V'_{i+1} \qquad (0 \leq i \leq D), 
\end{align}
where each of $V_{d+1} , \, V'_{D+1}$ is equal to $0$.
\item[\rm (iii)] The restrictions
\begin{align}
\label{eq:bp3}
[A', A]^{d-2i} |_{V_i} : V_i \rightarrow V_{d-i} \qquad (0 \leq i \leq d/2), \\
\label{eq:bp4}
[A, A']^{D-2i} |_{V'_i} : V'_i \rightarrow V'_{D-i} \qquad (0 \leq i \leq D/2), 
\end{align}
are bijections.
\end{enumerate}
\end{definition}

\begin{lemma}
\cite[Lemma 2.4]{Funk-Neubauer13}
\label{thm:dequalsdelta}
With reference to Definition \ref{def:bdpair}, we have $d=D$.
\end{lemma}

\begin{definition}
\rm
\label{def:bdtriple}
A {\it bidiagonal triple (BD triple) on $V$} is an ordered triple of linear transformations $A:V\to V$, $A' :V \to V$, and $A'' :V \to V$ that satisfy the following three conditions.
\begin{enumerate}
\item[\rm (i)] Each of $A,\, A', \, A''$ is diagonalizable.
\item[\rm (ii)] There exists an ordering $\lbrace V_i \rbrace _{i=0}^d$ (resp.~$\lbrace V'_i \rbrace _{i=0}^{D}$) (resp.~$\lbrace V''_i \rbrace _{i=0}^{\delta}$) of the eigenspaces of $A$ (resp.~$A'$) (resp. ~$A''$) with 
\begin{align}
\label{eq:bt1}
A' V_i \subseteq V_i + V_{i+1},  \qquad A'' V_i \subseteq V_{i-1} + V_i \qquad (0 \leq i \leq d), \\
\label{eq:bt2}
A'' V'_i \subseteq V'_i + V'_{i+1},  \qquad A V'_i \subseteq V'_{i-1} + V'_i \qquad (0 \leq i \leq D), \\
\label{eq:bt3}
AV''_i \subseteq V''_i + V''_{i+1},  \qquad A' V''_i \subseteq V''_{i-1} + V''_i \qquad (0 \leq i \leq \delta), 
\end{align}
where each of $V_{d+1}, \, V_{-1},\, V'_{D+1}, \, V'_{-1}, \, V''_{\delta + 1}, \, V''_{-1}$ is equal to $0$.
\item[\rm (iii)] The restrictions
\begin{align}
\label{eq:bt4}
[A', A]^{d-2i} |_{V_i} : V_i \rightarrow V_{d-i}, \qquad  [A'', A]^{d-2i} |_{V_{d-i}} : V_{d-i} \rightarrow V_i \qquad (0 \leq i \leq d/2), \\
\label{eq:bt5}
[A'', A']^{D-2i} |_{V'_i} : V'_i \rightarrow V'_{D-i}, \qquad  [A, A']^{D-2i} |_{V'_{D-i}} : V'_{D-i} \rightarrow V'_i \,\,\,\,\,\, (0 \leq i \leq D/2), \\
\label{eq:bt6}
[A, A'']^{\delta-2i} |_{V''_i} : V''_i \rightarrow V''_{\delta-i}, \qquad  [A', A'']^{\delta-2i} |_{V''_{\delta-i}} : V''_{\delta-i} \rightarrow V''_i  \qquad (0 \leq i \leq \delta/2),
\end{align}
are bijections.
\end{enumerate}
\end{definition}

\begin{lemma}
\label{thm:dDdelta}
With reference to Definition \ref{def:bdtriple}, we have $d=D=\delta$.
\end{lemma}

\noindent
{\it Proof:}  Let $A, \, A', \, A''$ denote a BD triple.  Then $A, \, A'$ is a BD pair.  So by Lemma \ref{thm:dequalsdelta} we have $d=D$.  Also, $A', \, A''$ is a BD pair.  So by Lemma \ref{thm:dequalsdelta} we have $D=\delta$.
\hfill $\Box$ \\

In view of Lemma \ref{thm:dDdelta}, for the remainder of this paper we use $d$ to index the eigenspaces of $A$, $A'$, and $A''$.  With reference to Definition \ref{def:bdpair} and Definition \ref{def:bdtriple} we call $V$ the {\it vector space underlying} $A, \, A'$ (resp. $A, \, A', \, A''$).  We say $A, \, A'$ (resp. $A, \, A', \, A''$) is {\it over} $\K$.  We call $d$ the {\it diameter} of $A, \, A'$ (resp. $A, \, A', \, A''$). \\

For the remainder of this paper we assume that the field $\K$ is algebraically closed and characteristic zero.  

\begin{definition}
\rm
\label{def:isom}
Let $A, \, A', \, A''$ and $B, \, B', \, B''$ denote BD triples over $\K$.  Let $V$ (resp.~$W$) denote the vector space underlying $A, \, A', \, A''$ (resp.~$B, \, B', \, B''$).  An {\it isomorphism of BD triples from $A, \, A', \, A''$ to $B, \, B', \, B''$} is a vector space isomorphism $\mu : V \rightarrow W$ such that $\mu A = B \mu$, $\mu A' = B' \mu$, and $\mu A'' = B'' \mu$.  We say $A, \, A', \, A''$ and $B, \, B', \, B''$ are {\it isomorphic} whenever there exists an isomorphism of BD triples from $A, \, A', \, A''$ to $B, \, B', \, B''$.
\end{definition}

The corresponding notion of isomorphism of BD pairs appears in \cite[Definition 2.11]{Funk-Neubauer13}. \\

Let $A, \, A', \, A''$ denote a BD triple of diameter $d$.  An ordering of the eigenspaces of $A$ (resp.~$A'$) (resp.~$A''$) is called {\it standard} whenever this ordering satisfies (\ref{eq:bt1}) (resp.~(\ref{eq:bt2})) (resp.~(\ref{eq:bt3})).  Let $\lbrace V_i \rbrace _{i=0}^d$ (resp.~$\lbrace V'_i \rbrace _{i=0}^d$) (resp.~$\lbrace V''_i \rbrace _{i=0}^d$) denote a standard ordering of the eigenspaces of $A$ (resp.~$A'$) (resp.~$A''$).  Then no other ordering of these eigenspaces is standard, and so the BD triple $A, \, A', \, A''$ uniquely determines these three standard orderings.  The corresponding notion of standard orderings of a BD pair appears in \cite[Section 2]{Funk-Neubauer13}.  For $0 \leq i \leq d$, let $\theta_i$ (resp.~$\theta'_i $) (resp.~$\theta''_i $) denote the eigenvalue of $A$ (resp.~$A'$) (resp.~$A''$) corresponding to $V_i$ (resp.~$V'_i$) (resp.~$V''_i$).  We call $\lbrace \theta_i \rbrace _{i=0}^d$ (resp.~$\lbrace \theta'_i \rbrace _{i=0}^d$) (resp.~$\lbrace \theta''_i \rbrace _{i=0}^d$) the {\it first} (resp.~{\it second}) (resp.~{\it third}) {\it eigenvalue sequence} of $A, \, A', \, A''$.  The corresponding notion of eigenvalue and dual eigenvalue sequence of a BD pair appears in \cite[Section 2]{Funk-Neubauer13}.  For the remainder of this section, we adopt the notation from this paragraph.    

\begin{lemma}
\label{thm:samedim}
For $0 \leq i \leq d$, the spaces $V_i$, $V_{d-i}$, $V'_i$, $V'_{d-i}$, $V''_i$, $V''_{d-i}$ all have the same dimension.
\end{lemma}

\noindent
{\it Proof:}  Since $A, \, A'$ is a BD pair then, by \cite[Lemma 2.7]{Funk-Neubauer13}, the spaces $V_i$, $V_{d-i}$, $V'_i$, $V'_{d-i}$ all have the same dimension.  Since $A', \, A''$ is a BD pair then, by \cite[Lemma 2.7]{Funk-Neubauer13}, the spaces $V'_i$, $V'_{d-i}$, $V''_i$, $V''_{d-i}$ all have the same dimension.
\hfill $\Box$ \\

\begin{definition}
\rm
\label{def:shape}
With reference to Lemma \ref{thm:samedim}, for $0 \leq i \leq d$, let $\rho_i$ denote the common dimension of  $V_i$, $V_{d-i}$, $V'_i$, $V'_{d-i}$, $V''_i$, $V''_{d-i}$.  We refer to the sequence $\lbrace \rho_i \rbrace _{i=0}^d$ as the {\it shape} of $A, \, A', \, A''$.
\end{definition}

The corresponding notion of the shape of a BD pair appears in \cite[Definition 2.8]{Funk-Neubauer13}.  

\begin{definition}
\rm
\label{def:parray}
The {\it parameter array} of $A, \, A', \, A''$ is the sequence \\ $(\lbrace \theta_i \rbrace _{i=0}^d; \lbrace \theta'_i \rbrace _{i=0}^d; \lbrace \theta''_i \rbrace _{i=0}^d; \lbrace \rho_i \rbrace _{i=0}^d)$, where $\lbrace \theta_i \rbrace _{i=0}^d$ (resp.~$\lbrace \theta'_i \rbrace _{i=0}^d$) (resp.~$\lbrace \theta''_i \rbrace _{i=0}^d$) is the first (resp. second) (resp. third) eigenvalue sequence of $A, \, A', \, A''$, and $\lbrace \rho_i \rbrace _{i=0}^d$ is the shape of $A, \, A', \, A''$.
\end{definition}

The corresponding notion of the parameter array of a BD pair appears in \cite[Definition 2.9]{Funk-Neubauer13}.  

\begin{lemma}
\label{thm:isompa}
Let $B, \, B', \, B''$ denote a BD triple over $\K$.  Then  $A, \, A', \, A''$ and $B, \, B', \, B''$ are isomorphic if and only if the parameter array of $A, \, A', \, A''$ equals the parameter array of $B, \, B', \, B''$.
\end{lemma}

We postpone the proof of Lemma \ref{thm:isompa} until Section 7. 

\begin{lemma}
\label{thm:recur}
Suppose that $d \geq 2$.  Then the expressions 
\begin{eqnarray}
\label{recur} 
\frac{\theta_{i+1} - \theta_i}{\theta_i - \theta_{i-1}}, \qquad \frac{\theta'_{i+1} - \theta'_i}{\theta'_i - \theta'_{i-1}}, \qquad \frac{\theta''_{i+1} - \theta''_i}{\theta''_i - \theta''_{i-1}}
\end{eqnarray}
are equal and independent of $i$ for $1 \leq i \leq d-1$.
\end{lemma}

\noindent
{\it Proof:}  Observe that $A, A'$ is a BD pair with eigenvalue (resp.~dual eigenvalue) sequence $\lbrace \theta_i \rbrace _{i=0}^d$ (resp.~$\lbrace \theta'_{d-i} \rbrace _{i=0}^d$).  So, by \cite[Theorem 5.1]{Funk-Neubauer13}, 
\begin{eqnarray*}
\frac{\theta_{i+1} - \theta_i}{\theta_i - \theta_{i-1}}, \qquad \frac{\theta'_{i+1} - \theta'_i}{\theta'_i - \theta'_{i-1}}
\end{eqnarray*}
are equal and independent of $i$ for $1 \leq i \leq d-1$.
Also, $A', A''$ is a BD pair with eigenvalue (resp. dual eigenvalue) sequence $\lbrace \theta'_i \rbrace _{i=0}^d$ (resp.~$\lbrace \theta''_{d-i} \rbrace _{i=0}^d$).  So, by \cite[Theorem 5.1]{Funk-Neubauer13}, 
\begin{eqnarray*}
\frac{\theta'_{i+1} - \theta'_i}{\theta'_i - \theta'_{i-1}}, \qquad \frac{\theta''_{i+1} - \theta''_i}{\theta''_i - \theta''_{i-1}}
\end{eqnarray*}
are equal and independent of $i$ for $1 \leq i \leq d-1$.
\hfill $\Box$ \\

\begin{definition}
\rm
\label{def:base}
We refer to the common value of (\ref{recur}) as the {\it base} of $A, \, A', \, A''$.  By Lemma \ref{thm:recur}, the base is defined for $d \geq 2$.  For $d \leq 1$, we adopt the convention that the base is $1$.  
\end{definition}

Since $\lbrace \theta_i \rbrace _{i=0}^d$ is a list of {\it distinct} scalars then the base of $A, \, A', \, A''$ is always nonzero.  The corresponding notion of the base of a BD pair appears in \cite[Definition 5.5]{Funk-Neubauer13}.  See also \cite[Lemma 9.1]{Funk-Neubauer13}.  

\begin{note}
\label{bnot1}
\rm
Let $b$ denote the base of $A, \, A', \, A''$.  Whenever $b \neq 1$, we let $q \in \K$ denote a root of the polynomial $\lambda ^2 - b^{-1} \in \K[\lambda]$.  The scalar $q$ exists since $\K$ is algebraically closed.  Observe $q \neq 0$, and so $b = q^{-2}$.  By construction $b$ uniquely determines $q$ up to sign.
\end{note}

The following notions of affine equivalence will be used in stating some of our main results (see Theorem \ref{thm:unique} and Theorem \ref{thm:reduce}).  

\begin{definition}
\rm
\label{def:ae}
Let $\lbrace \sigma_i \rbrace _{i=0}^d$ and $\lbrace \tau_i \rbrace _{i=0}^d$ each denote a sequence of scalars taken from $\K$.  Let $X:V\to V$ and $Y:V \to V$ denote linear transformations.  Let $B, \, B', \, B''$ denote a BD triple on $V$.
\begin{enumerate}
\item[\rm (i)] We say $\lbrace \sigma_i \rbrace _{i=0}^d$ is {\it affine equivalent} to $\lbrace \tau_i \rbrace _{i=0}^d$, denoted $\lbrace \sigma_i \rbrace _{i=0}^d \sim \lbrace \tau_i \rbrace _{i=0}^d$, whenever there exist $r, s$ in $\K$ with $r \neq 0$ such that $\sigma_i = r \tau_i + s$ ($0 \leq i \leq d$).  
\item[\rm (ii)] We say $X$ is {\it affine equivalent} to $Y$, denoted $X \sim Y$, whenever there exist $r, s$ in $\K$ with $r \neq 0$ such that $X = r Y + sI$.  
\item[\rm (iii)] We say $A, \, A', \, A''$ is {\it affine equivalent} to $B, \, B', \, B''$ whenever $A \sim B$, $A' \sim B'$, and $A'' \sim B''$.  
\end{enumerate}
\end{definition}

\section{The Algebras $\SL$ and $\uq$ }

First we recall the Lie algebra $\SL$ and its finite-dimensional modules.

\begin{definition}
\rm
\label{def:usl2}
Let $\SL$ denote the Lie algebra over $\K$ that has a basis $h, \, e, \, f$ and Lie bracket
\begin{eqnarray*}
[h,e] = 2e, \qquad
[h,f] = -2f, \qquad
[e,f] = h.
\end{eqnarray*}
\end{definition}

\begin{theorem}
\cite[Lemma 3.2]{HarTer07}
\label{thm:ueq}
The Lie algebra $\SL$ is isomorphic to the Lie algebra over $\K$ that has basis $X,\, Y, \, Z$ and Lie bracket
\begin{eqnarray*}
[X,Y] = 2X + 2Y, \qquad
[Y,Z] = 2Y + 2Z, \qquad
[Z,X] =2Z + 2X.
\end{eqnarray*}
An isomorphism with the presentation in Definition \ref{def:usl2} is given by:
\begin{eqnarray*}
X &\rightarrow& 2e - h, \\
Y &\rightarrow& -2f - h, \\
Z &\rightarrow& h.
\end{eqnarray*}
The inverse of this isomorphism is given by:
\begin{eqnarray*}
e &\rightarrow& (X + Z)/2, \\
f &\rightarrow& -(Y + Z)/2, \\
h &\rightarrow& Z.
\end{eqnarray*}
\end{theorem}

\begin{definition}
\rm
\label{def:equitpair}
Given an ordered triple $X, \, Y, \, Z$ of elements in $\SL$, we call this triple {\it equitable} whenever $X, \, Y, \, Z$ satisfy the relations from Theorem \ref{thm:ueq}.
\end{definition}

The following two lemmas give a description of all finite-dimensional $\SL$-modules.

\begin{lemma}
\cite[Theorem 6.3]{Humph72}
\label{thm:compred}
Each finite-dimensional $\SL$-module $V$ is completely reducible; this means that $V$ is a direct sum of irreducible $\SL$-modules.
\end{lemma}

The finite-dimensional irreducible $\SL$-modules are described as follows.

\begin{lemma}
\cite[Theorem 7.2]{Humph72}
\label{thm:sl2mods}
There exists a family of finite-dimensional irreducible $\SL$-modules
\begin{eqnarray*}
V(d), \qquad  d=0,1,2,\ldots
\end{eqnarray*}
with the following properties:  $V(d)$ has a basis $\lbrace v_i \rbrace _{i=0}^d$ such that $h.v_i = (d-2i) v_i$ for $0 \leq i \leq d$, $f.v_i = (i+1) v_{i+1}$ for $0 \leq i \leq d$, where $v_{d+1}=0$, and $e.v_i = (d-i+1) v_{i-1}$ for $0 \leq i \leq d$, where $v_{-1}=0$.  Moreover, every finite-dimensional irreducible $\SL$-module is isomorphic to exactly one of the modules $V(d)$.
\end{lemma}

\begin{lemma}
\cite[Lemma 3.7]{Funk-Neubauer13}
\label{thm:2pieces}
Let $V$ denote an $\SL$-module with finite positive dimension (not necessarily irreducible).  Define 
\begin{eqnarray*}
V_{\rm{even}} &:=& \operatorname{span}\{\,v \in V\, |\, h.v = i \, v, \,\, i \in \Z, \, i \, even \,\}, \\
V_{\rm{odd}} &:=& \operatorname{span}\{\,v \in V\, |\, h.v = i \, v, \,\, i \in \Z, \, i \, odd \,\}.
\end{eqnarray*}
Then $V_{\rm{even}}$ and $V_{\rm{odd}}$ are $\SL$-modules, and $V = V_{\rm{even}} + V_{\rm{odd}}$ (direct sum).
\end{lemma}

The following definition will be used in stating two of our main results (see Theorem \ref{thm:uqsl2bidiag} and Theorem \ref{thm:reducedtomod}).

\begin{definition}
\rm
\label{def:pure}
Let $V$ denote an $\SL$-module with finite positive dimension.  With reference to Lemma \ref{thm:2pieces}, we say $V$ is {\it segregated} whenever $V=V_{\rm{even}}$ or $V=V_{\rm{odd}}$.
\end{definition}

We now recall the quantum algebra $\uq$ and its finite-dimensional modules.  Let $q$ denote a nonzero scalar in $\K$ which is not a root of unity.  

\begin{definition}
\rm
\label{def:uqsl2}
Let $\uq$ denote the unital associative $\K$-algebra with generators \\ $k, \, k^{-1}, \, e, \, f$ and the following relations:
\begin{eqnarray*}
kk^{-1} &=& k^{-1}k = 1, \\
ke &=& q^2 ek, \\
kf &=& q^{-2}fk, \\
ef - fe &=& \frac{k - k^{-1}}{q - q^{-1}}.
\end{eqnarray*}
\end{definition}

\begin{theorem}
\cite[Theorem 2.1]{ItoTerWang06}
\label{thm:uqeq}
The algebra $\uq$ is isomorphic to the unital associative $\K$-algebra with generators $x, \, x^{-1}, \, y, \, z$ and the following relations:
\begin{eqnarray*}
xx^{-1} = x^{-1}x &=& 1, \\
\frac{qxy - q^{-1}yx}{q - q^{-1}} &=& 1, \\
\frac{qyz - q^{-1}zy}{q - q^{-1}} &=& 1, \\
\frac{qzx - q^{-1}xz}{q - q^{-1}} &=& 1.
\end{eqnarray*}
An isomorphism with the presentation in Definition \ref{def:uqsl2} is given by:
\begin{eqnarray*}
x^{\pm1} &\rightarrow& k^{\pm1}, \\
y &\rightarrow& k^{-1} + f(q-q^{-1}), \\
z &\rightarrow& k^{-1} - k^{-1}eq(q-q^{-1}).
\end{eqnarray*}
The inverse of this isomorphism is given by:
\begin{eqnarray*}
k^{\pm1} &\rightarrow& x^{\pm1}, \\
f &\rightarrow& (y - x^{-1})(q-q^{-1}), \\
e &\rightarrow& (1 - xz)q^{-1}(q-q^{-1}).
\end{eqnarray*}
\end{theorem}

\begin{definition}
\rm
\label{def:equitpair2}
Given an ordered triple $x, \, y, \, z$ of elements in $\uq$, we call this triple {\it equitable} whenever $x, \, y, \, z$ satisfy the relations from Theorem \ref{thm:uqeq}.
\end{definition}

The following two lemmas give a description of all finite-dimensional $\uq$-modules.

\begin{lemma}
\cite[Theorems 2.3, 2.9]{Jantzen96}
\label{thm:uqcompred}
Each finite-dimensional $\uq$-module $V$ is completely reducible; this means that $V$ is a direct sum of irreducible $\uq$-modules.
\end{lemma}

For each nonnegative integer $n$, define $[n] := (q^{n} - q^{-n})/(q-q^{-1})$.  The finite-dimensional irreducible $\uq$-modules are described as follows.  

\begin{lemma}
\cite[Theorem 2.6]{Jantzen96}
\label{thm:uq2mods}
There exists a family of finite-dimensional irreducible $\uq$-modules
\begin{eqnarray*}
V(d, \epsilon), \qquad \epsilon \in \{1, -1 \}, \qquad  d=0,1,2,\ldots
\end{eqnarray*}
with the following properties:  $V(d, \epsilon)$ has a basis $\lbrace v_i \rbrace _{i=0}^d$ such that $k.v_i = \epsilon q^{d-2i} v_i$ for $0 \leq i \leq d$, $f.v_i = [i+1] v_{i+1}$ for $0 \leq i \leq d$, where $v_{d+1}=0$, and $e.v_i = \epsilon [d-i+1] v_{i-1}$ for $0 \leq i \leq d$, where $v_{-1}=0$.  Moreover, every finite-dimensional irreducible $\uq$-module is isomorphic to exactly one of the modules $V(d, \epsilon)$.
\end{lemma}

\begin{lemma}
\cite[Lemma 4.7]{Funk-Neubauer13}
\label{thm:4pieces}
Let $V$ denote a $\uq$-module with finite positive dimension (not necessarily irreducible).  For $\epsilon \in \{1, -1 \}$, define
\begin{eqnarray*}
V_{\rm{even}}^{\epsilon} &:=& \operatorname{span}\{\,v \in V\, |\, k.v = \epsilon \, q^i \, v, \,\, i \in \Z, \, i \, even \,\}, \\
V_{\rm{odd}}^{\epsilon} &:=& \operatorname{span}\{\,v \in V\, |\, k.v = \epsilon \, q^i \, v, \,\, i \in \Z, \, i \, odd \,\}.
\end{eqnarray*}
Then $V_{\rm{even}}^{1}, \, V_{\rm{even}}^{-1}, \, V_{\rm{odd}}^{1}, \, V_{\rm{odd}}^{-1}$ are $\uq$-modules, and $V = V_{\rm{even}}^{1} + V_{\rm{even}}^{-1} + V_{\rm{odd}}^{1} + V_{\rm{odd}}^{-1}$ (direct sum).
\end{lemma}

The following definition will be used in stating two of our main results (see Theorem \ref{thm:uqsl2bidiag} and Theorem \ref{thm:reducedtomod}).

\begin{definition}
\rm
\label{def:uqpure}
Let $V$ denote a $\uq$-module with finite positive dimension.  With reference to Lemma \ref{thm:4pieces}, we say  $V$ is {\it segregated} whenever $V=V_{\rm{even}}^1$ or $V=V_{\rm{odd}}^1$.
\end{definition}

\section{The Main Theorems}

The seven theorems in this section make up the main conclusions of the paper.  The following theorem states that every BD pair can be extended to a BD triple. 

\begin{theorem}
\label{thm:extend}
Let $A, \, A'$ denote a BD pair on $V$ with base $b$ and parameter array \\ $(\lbrace \theta_i \rbrace _{i=0}^d; \lbrace \theta'_i \rbrace _{i=0}^d; \lbrace \rho_i \rbrace _{i=0}^d)$.  Then there exists a linear transformation $A'': V \rightarrow V$ such that $A, \, A', \, A''$ is a BD triple on $V$ of base $b$.  Moreover,  $\lbrace \theta_i \rbrace _{i=0}^d$ (resp. $\lbrace \theta'_{d-i} \rbrace _{i=0}^d$) is the first (resp. second) eigenvalue sequence of $A, \, A', \, A''$, and $\lbrace \rho_i \rbrace _{i=0}^d$ is the shape of $A, \, A', \, A''$.
\end{theorem}

We refer to Theorem \ref{thm:extend} as the extension theorem.  In the proof of Theorem \ref{thm:extend} we explicitly construct the linear transformation $A''$ (see Section 6).  The following theorem implies that the linear transformation $A''$ from Theorem \ref{thm:extend} is uniquely determined up to affine equivalence.  See Corollary \ref{thm:extuniq}. 

\begin{theorem}
\label{thm:unique}
Let $A, \, A', \, A''$ and $B, \, B', \, B''$ denote BD triples on $V$.  Assume \\ $A \sim B$ and $A' \sim B'$.  Then $A'' \sim B''$.  
\end{theorem}

We refer to Theorem \ref{thm:unique} as the uniqueness theorem.  The following theorem provides a classification of BD triples up to isomorphism.

\begin{theorem}
\label{thm:class}
Let
\begin{eqnarray}
\label{parameterarray}
(\lbrace \theta_i \rbrace _{i=0}^d; \lbrace \theta'_i \rbrace _{i=0}^d; \lbrace \theta''_i \rbrace _{i=0}^d; \lbrace \rho_i \rbrace _{i=0}^d)
\end{eqnarray}
denote a sequence of scalars taken from $\K$.  Then there exists a BD triple $A, \, A', \, A''$ over $\K$ with parameter array (\ref{parameterarray}) if and only if the following {\rm(i)}--{\rm(v)} hold.
\begin{enumerate}
\item[\rm (i)] $\theta_i \neq \theta_j$, $\theta'_i \neq \theta'_j$, $\theta''_i \neq \theta''_j$ for $0 \leq i,j \leq d$ and $i \neq j$.
\item[\rm (ii)] The expressions
\begin{eqnarray*}
\frac{\theta_{i+1} - \theta_i}{\theta_i - \theta_{i-1}}, \qquad \frac{\theta'_{i+1} - \theta'_i}{\theta'_i - \theta'_{i-1}}, \qquad \frac{\theta''_{i+1} - \theta''_i}{\theta''_i - \theta''_{i-1}},
\end{eqnarray*}
are equal and independent of $i$ for $1 \leq i \leq d-1$.
\item[\rm (iii)] $\rho_i$ is a positive integer for $0 \leq i \leq d$.
\item[\rm (iv)] $\rho_i = \rho_{d-i}$ for $0 \leq i \leq d$.
\item[\rm (v)] $\rho_{i} \leq \rho_{i+1}$ for $0 \leq i < d/2$.
\end{enumerate}
Suppose that {\rm(i)--(v)} hold.  Then $A, \, A', \, A''$ is unique up to isomorphism of BD triples.
\end{theorem}

We refer to Theorem \ref{thm:class} as the classification theorem.  With reference to Theorem \ref{thm:class}, for $d=0$ we regard conditions (i), (ii), (v) as holding since they are vacuously true.  Similarly, for $d=1$ we regard condition (ii) as holding since it is vacuously true.  It is not clear from the statement of Theorem \ref{thm:class} how BD triples are related to finite-dimensional modules for $\SL$ and $\uq$.  Also, it is not clear from Theorem \ref{thm:class} how the BD triples in each isomorphism class are constructed.  The next four theorems address these issues. 

\begin{theorem}
\label{thm:triplerelations}
Let $A, \, A', \, A''$ denote a BD triple of base $b$.  Then there exists a sequence of scalars $\alpha, \, \alpha', \, \alpha'', \, \gamma_1, \, \gamma_2, \, \gamma_3$ in $\K$ such that
\begin{eqnarray}
\label{eq:triple1}
A A' - b A' A - \alpha' A - \alpha A' - \gamma_1 I = 0, \\
\label{eq:triple2}
A' A'' - b A'' A' - \alpha'' A' - \alpha' A'' - \gamma_2 I = 0, \\
\label{eq:triple3}
A'' A - b A A'' - \alpha A'' - \alpha'' A - \gamma_3 I = 0.
\end{eqnarray}
This sequence of scalars is uniquely determined by the triple $A, \, A', \, A''$ provided $d \geq 2$.
\end{theorem}

We refer to (\ref{eq:triple1})--(\ref{eq:triple3}) as the fundamental bidiagonal relations.  See Lemma \ref{thm:6recurs} for a description of how the scalars $b, \, \alpha, \, \alpha', \, \alpha'', \, \gamma_1, \, \gamma_2, \, \gamma_3$ are related to the eigenvalues of $A, \, A', \, A''$.  The following definition will be used to state the next three theorems.

\begin{definition}
\label{def:reduced}
\rm
Let $A, \, A', \, A''$ denote a BD triple of diameter $d$ and base $b$.  We say $A, \, A', \, A''$ is {\it reduced} if either (i) or (ii) holds.
\begin{enumerate}
\item[\rm (i)] $b=1$ and the first, second, and third eigenvalue sequences of $A, \, A', \, A''$ are each $\lbrace 2i-d \rbrace _{i=0}^d$.
\item[\rm (ii)] $b \neq 1$ and the first, second, and third eigenvalue sequences of $A, \, A', \, A''$ are each $\lbrace q^{d-2i} \rbrace _{i=0}^d$.
\end{enumerate}
\end{definition}

\begin{theorem}
\label{thm:reduce}
Every BD triple is affine equivalent to a reduced BD triple.
\end{theorem}

We refer to Theorem \ref{thm:reduce} as the reducibility theorem. 

\begin{note}
\label{reducedcases}
\rm
We make the following observation in order to motivate the next two theorems.  Let $A, \, A', \, A''$ denote a reduced BD triple of base $b$.  It is shown in Corollary \ref{thm:reducedrels} that for $b=1$, (\ref{eq:triple1})--(\ref{eq:triple3}) take the form
\begin{eqnarray*}
A A' - A' A - 2A - 2 A' = 0, \\
A' A'' - A'' A' - 2A' - 2 A'' = 0, \\
A'' A - A A'' - 2A'' - 2 A = 0.
\end{eqnarray*}
It is also shown in Corollary \ref{thm:reducedrels} that for $b \neq 1$, (\ref{eq:triple1})--(\ref{eq:triple3}) take the form
\begin{eqnarray*}
q A A' - q^{-1} A' A - (q - q^{-1}) I =  0, \\
q A' A'' - q^{-1} A'' A' - (q - q^{-1}) I =  0, \\
q A'' A - q^{-1} A A'' - (q - q^{-1}) I =  0.
\end{eqnarray*}
\end{note}

The following two theorems provide a correspondence between reduced BD triples and segregated $\SL$/$\uq$-modules.
\begin{theorem}
\label{thm:uqsl2bidiag}
Let $V$ denote a segregated $\SL$-module (resp.~segregated $\uq$-module).  Then each equitable triple in $\SL$ (resp.~$\uq$) acts on $V$ as a reduced BD triple of base $1$ (resp.~base not equal to $1$).
\end{theorem}

The next theorem can be thought of as a converse of Theorem \ref{thm:uqsl2bidiag}.

\begin{theorem}
\label{thm:reducedtomod}
Let $A, \, A', \, A''$ denote a reduced BD triple on $V$, and let $b$ denote its base.  Then the following {\rm(i)},{\rm(ii)} hold.
\begin{enumerate}
\item[\rm (i)] Suppose that $b = 1$.  Let $X, \, Y, \, Z$ denote an equitable triple in $\SL$.  Then there exists an $\SL$-module structure on $V$ such that $(X - A)V=0$, $(Y - A')V=0$, and $(Z - A'')V=0$.  This $\SL$-module structure on $V$ is segregated.
\item[\rm (ii)] Suppose that $b \neq 1$.  Let $x, \, y, \, z$ denote an equitable triple in $\uq$.  Then there exists a $\uq$-module structure on $V$ such that $(x - A)V=0$, $(y - A')V=0$, and $(z - A'')V=0$.  This $\uq$-module structure on $V$ is segregated.
\end{enumerate}
\end{theorem}

We refer to Theorems \ref{thm:uqsl2bidiag} and \ref{thm:reducedtomod} as the correspondence theorems.  A version of Theorem \ref{thm:reducedtomod} in which $A, \, A', \, A''$ is assumed to have shape $(1, 1, 1, \ldots, 1)$ appears in \cite[Theorem 18.3]{Ter13}.  Theorem \ref{thm:reduce} and Theorem \ref{thm:reducedtomod} together show that every BD triple is affine equivalent to a BD triple of the type constructed in Theorem \ref{thm:uqsl2bidiag}.

\section{Preliminaries}

In this section we develop the tools needed to prove our main results.  

\begin{definition}
\label{def:brec}
\rm
Let $b \in \K$.  Suppose that $d \geq 2$ and let $\lbrace \sigma_i \rbrace _{i=0}^d$ denote a sequence of distinct scalars taken from $\K$.  This sequence is called {\it $b$-recurrent} whenever
\begin{eqnarray*}
\frac{\sigma_{i+1} - \sigma_i}{\sigma_i - \sigma_{i-1}} = b \qquad (1 \leq i \leq d-1).
\end{eqnarray*}
\end{definition}

\begin{lemma}
\label{thm:brecuraffine}
Suppose that $d \geq 2$ and let $\lbrace \sigma_i \rbrace _{i=0}^d$ and $\lbrace \tau_i \rbrace _{i=0}^d$ each denote a sequence of distinct scalars taken from $\K$.  Then the following {\rm(i)},{\rm(ii)} hold.  
\begin{enumerate}
\item[\rm (i)] Assume $\lbrace \sigma_i \rbrace _{i=0}^d$ and $\lbrace \tau_i \rbrace _{i=0}^d$ are each $b$-recurrent.  Then $\lbrace \sigma_i \rbrace _{i=0}^d \sim \lbrace \tau_i \rbrace _{i=0}^d$.   
\item[\rm (ii)] Assume $\lbrace \sigma_i \rbrace _{i=0}^d$ is $b$-recurrent and $\lbrace \sigma_i \rbrace _{i=0}^d \sim \lbrace \tau_i \rbrace _{i=0}^d$.  Then $\lbrace \tau_i \rbrace _{i=0}^d$ is $b$-recurrent.
\end{enumerate}
\end{lemma}

\noindent
{\it Proof:}  (i) Solving the recurrence from Definition \ref{def:brec} we find that there exist $c_1, c_2$ in $\K$ with $c_2 \neq 0$ such that for $b=1$ (resp.~$b \neq 1$), $\sigma_i = c_1 + c_2 i$ (resp.~$\sigma_i = c_1 + c_2 b^i$) ($0 \leq i \leq d$).   Similarly, there exist $c_3, c_4$ in $\K$ with $c_4 \neq 0$ such that for $b=1$ (resp.~$b \neq 1$), $\tau_i = c_3 + c_4 i$ (resp.~$\tau_i = c_3 + c_4 b^i$) ($0 \leq i \leq d$).  Hence, $\sigma_i = c_2 c_4 ^{-1} \tau_i + c_1 - c_2 c_3 c_4 ^{-1}$ ($0 \leq i \leq d$), and the result follows.   \\
(ii) By assumption there exist $r, s$ in $\K$ with $r \neq 0$ such that $\sigma_i = r \tau_i + s$ ($0 \leq i \leq d$).  From this and since $\lbrace \sigma_i \rbrace _{i=0}^d$ is $b$-recurrent we have 
\begin{eqnarray*}
b = \frac{\sigma_{i+1} - \sigma_i}{\sigma_i - \sigma_{i-1}} = \frac{r \tau_{i+1} + s - (r \tau_i + s)}{r \tau_i + s - (r \tau_{i-1} + s)} = \frac{\tau_{i+1} - \tau_i}{\tau_i - \tau_{i-1}} \qquad (1 \leq i \leq d-1).
\end{eqnarray*} 
Thus, $\lbrace \tau_i \rbrace _{i=0}^d$ is $b$-recurrent.
\hfill $\Box$ \\

\begin{lemma}
\label{thm:d<2}
Suppose that $d < 2$ and let $\lbrace \sigma_i \rbrace _{i=0}^d$ and $\lbrace \tau_i \rbrace _{i=0}^d$ each denote a sequence of distinct scalars taken from $\K$.  Then $\lbrace \sigma_i \rbrace _{i=0}^d \sim \lbrace \tau_i \rbrace _{i=0}^d$.   
\end{lemma}

\noindent
{\it Proof:}  First assume that $d=1$.  Let $r := (\sigma_0 - \sigma_1)(\tau_0 - \tau_1)^{-1}$ and $s := \sigma_0 - r \tau_0$.  Observe $r \neq 0$ since $\sigma_0 \neq \sigma_1$.  We now have $\sigma_i = r \tau_i + s$ for $0 \leq i \leq 1$, and so $\lbrace \sigma_i \rbrace _{i=0}^1 \sim \lbrace \tau_i \rbrace _{i=0}^1$.  Now assume that $d=0$.  Let $r := 1$ and $s := \sigma_0 - \tau_0$.  Hence, $\sigma_0 = r \tau_0 + s$, and so $\lbrace \sigma_i \rbrace _{i=0}^0 \sim \lbrace \tau_i \rbrace _{i=0}^0$.
\hfill $\Box$ \\

\begin{definition}
\label{def:decomp}
\rm
A {\it decomposition} of $V$ is a sequence  $\lbrace U_i \rbrace _{i=0}^d$ consisting of nonzero subspaces of $V$ such that $V=\sum_{i=0}^{d} U_i$ (direct sum).  For any decomposition of $V$ we adopt the convention that $U_{-1} := 0$ and $U_{d+1} := 0$.  
\end{definition}

\begin{lemma}
\label{thm:mapseq}
Let $\lbrace U_i \rbrace _{i=0}^d$ denote a decomposition of $V$.  Let $\lbrace \sigma_i \rbrace _{i=0}^d$ and $\lbrace \tau_i \rbrace _{i=0}^d$ each denote a sequence of distinct scalars taken from $\K$.  Let $X :V \to V$ (resp.~$Y: V \to V$) denote the linear transformation such that for $0 \leq i \leq d$, $U_i$ is an eigenspace for $X$ (resp.~$Y$) with eigenvalue $\sigma_i$ (resp.~$\tau_i$).  Then $\lbrace \sigma_i \rbrace _{i=0}^d \sim \lbrace \tau_i \rbrace _{i=0}^d$ if and only if $X \sim Y$.   
\end{lemma}

\noindent
{\it Proof:}  $(\Longrightarrow)$:  Let $v \in V$.  Since $\lbrace U_i \rbrace _{i=0}^d$ is a decomposition of $V$ there exist $u_i \in U_i$ such that $v = \sum_{i=0}^d u_i$.  By assumption there exist $r, s$ in $\K$ with $r \neq 0$ such that $\sigma_i = r \tau_i + s$ $(0 \leq i \leq d)$.  From this we have  
\begin{align*}
Xv &= \sum_{i=0}^d \sigma_i u_i  \\
&= r \sum_{i=0}^d \tau_i u_i + s \sum_{i=0}^d u_i \\
&=  r \sum_{i=0}^d Y u_i + s \sum_{i=0}^d u_i \\
&= (rY + sI)v.
\end{align*} 
Thus, $X = rY + sI$, and so $X \sim Y$. \\
$(\Longleftarrow)$:  Similar to $(\Longrightarrow)$.
\hfill $\Box$ \\
 
 We now develop some more properties of BD triples.  Throughout the remainder of this paper we will refer to the following assumption.
 
\begin{assumption}
\label{assump}
\rm
Let $A, \, A', \, A''$ denote a BD triple on $V$ of diameter $d$ and base $b$.  Let $\lbrace V_i \rbrace _{i=0}^d$ (resp.~$\lbrace V'_i \rbrace _{i=0}^d$) (resp.~$\lbrace V''_i \rbrace _{i=0}^d$) denote the standard ordering of the eigenspaces of $A$ (resp.~$A'$) (resp.~$A''$).  Let $(\lbrace \theta_i \rbrace _{i=0}^d; \lbrace \theta'_i \rbrace _{i=0}^d; \lbrace \theta''_i \rbrace _{i=0}^d; \lbrace \rho_i \rbrace _{i=0}^d)$ denote the parameter array of $A, \, A', \, A''$.  
\end{assumption}

Referring to Definition \ref{def:decomp} and Assumption \ref{assump}, the sequences  $\lbrace V_i \rbrace _{i=0}^d$, $\lbrace V'_i \rbrace _{i=0}^d$, $\lbrace V''_i \rbrace _{i=0}^d$ are each decompositions of $V$.   

\begin{lemma}
\label{thm:affineshift}
Let $r, \, s, \, t, \, u, \, v, \, w$ denote scalars in $\K$ with $r, \, t, \, v$ each nonzero.  With reference to Assumption \ref{assump}, $rA + sI , \, tA' + uI, \, vA'' + wI$ is a BD triple on $V$ of base $b$ and with parameter array $(\lbrace r \theta_i + s \rbrace _{i=0}^d; \lbrace t \theta'_i  + u \rbrace _{i=0}^d; \lbrace v \theta''_i  + w \rbrace _{i=0}^d; \lbrace \rho_i \rbrace _{i=0}^d)$.
\end{lemma}

\noindent
{\it Proof:}  Imitate the proof of  \cite[Lemma 2.10]{Funk-Neubauer13}. 
\hfill $\Box$ \\

\begin{lemma}
\label{thm:triple3pairs}
With reference to Assumption \ref{assump}, the following {\rm(i)}--{\rm(iii)} hold.
\begin{enumerate}
\item[\rm (i)] $A, \, A'$ is a BD pair with parameter array $(\lbrace \theta_i \rbrace _{i=0}^d; \lbrace \theta'_{d-i} \rbrace _{i=0}^d; \lbrace \rho_i \rbrace _{i=0}^d)$.
\item[\rm (ii)] $A', \, A''$ is a BD pair with parameter array $(\lbrace \theta'_i \rbrace _{i=0}^d; \lbrace \theta''_{d-i} \rbrace _{i=0}^d; \lbrace \rho_i \rbrace _{i=0}^d)$.
\item[\rm (ii)] $A'', \, A$ is a BD pair with parameter array $(\lbrace \theta''_i \rbrace _{i=0}^d; \lbrace \theta_{d-i} \rbrace _{i=0}^d; \lbrace \rho_i \rbrace _{i=0}^d)$.
\end{enumerate}
\end{lemma}

\noindent
{\it Proof:}  Immediate from Definitions \ref{def:bdpair}, \ref{def:bdtriple}, \ref{def:parray} and \cite[Definition 2.9]{Funk-Neubauer13}.  
\hfill $\Box$ \\

\begin{lemma}
\label{thm:3cornersums}
With reference to Assumption \ref{assump}, the following {\rm(i)}--{\rm(iii)} hold.
\begin{enumerate}
\item[\rm (i)] $V_i + \cdots + V_d  = V'_0 + \cdots + V'_{d-i} \qquad (0 \leq i \leq d)$.
\item[\rm (ii)] $V'_i + \cdots + V'_d = V''_0 + \cdots + V''_{d-i} \qquad (0 \leq i \leq d)$. 
\item[\rm (iii)] $V''_i + \cdots + V''_d  = V_0 + \cdots + V_{d-i} \qquad (0 \leq i \leq d)$.
\end{enumerate}
\end{lemma}

\noindent
{\it Proof:}  (i) By Lemma \ref{thm:triple3pairs}(i), $A, \, A'$ is a BD pair and $\lbrace V_i \rbrace _{i=0}^d$ (resp.~$\lbrace V'_{d-i} \rbrace _{i=0}^d$) is an ordering of the eigenspaces of $A$ (resp.~$A'$) satisfying (\ref{eq:bp1}) (resp.~(\ref{eq:bp2})).  Combining this with \cite[Lemma 6.6]{Funk-Neubauer13} we obtain the desired result. \\
(ii), (iii)  Similar to (i). 
\hfill $\Box$ \\
 
\begin{lemma}
\label{thm:intersect}
With reference to Assumption \ref{assump}, the following {\rm(i)}--{\rm(iii)} hold.
\begin{enumerate}
\item[\rm (i)] $V_i = (V'_0 + \cdots + V'_{d-i}) \cap (V''_{d-i} + \cdots + V''_d) \qquad (0 \leq i \leq d)$.
\item[\rm (ii)] $V'_i = (V''_0 + \cdots + V''_{d-i}) \cap (V_{d-i} + \cdots + V_d) \qquad (0 \leq i \leq d)$.
\item[\rm (iii)] $V''_i = (V_0 + \cdots + V_{d-i}) \cap (V'_{d-i} + \cdots + V'_d) \qquad (0 \leq i \leq d)$.
\end{enumerate}
\end{lemma}
 
 \noindent
{\it Proof:}  (i) Immediate from Lemma \ref{thm:3cornersums}(i),(iii) and the fact that $\sum_{i=0}^d V_i$ is a direct sum. \\
(ii), (iii) Similar to (i).
\hfill $\Box$ \\
 
For the remainder of this section we let $W$ denote a vector space over $\K$ with finite positive dimension and $\mu : V \to W$ denote a vector space isomorphism. 
 
\begin{lemma}
\label{thm:muj}
Let $X: V \to V$ and $Y: V \to V$ denote linear transformations.  Let $j$ denote a nonnegative integer.  Then $[\mu X \mu^{-1}, \mu Y \mu^{-1}]^j = \mu [X, Y]^j \mu^{-1}$. 
\end{lemma}

\noindent
{\it Proof:}  The proof is by induction on $j$.  The result is trivially true for $j=0$.  Now assume that $j \geq 1$.  By induction we have that 
\begin{align*}
[\mu X \mu^{-1}, \mu Y \mu^{-1}]^{j+1} &=  [\mu X \mu^{-1}, \mu Y \mu^{-1}] \mu [X, Y]^j \mu^{-1}  \\
&= (\mu X Y \mu^{-1} - \mu Y X \mu^{-1}) \mu [X, Y]^j \mu^{-1} \\
&= \mu (X Y - Y X) [X, Y]^j \mu^{-1}   \\
&= \mu [X, Y]^{j+1} \mu^{-1},
\end{align*}     
and the result follows.
\hfill $\Box$ \\

\begin{lemma}
\label{thm:muiso}
With reference to Assumption \ref{assump}, $\mu A \mu^{-1}, \, \mu A' \mu^{-1}, \, \mu A'' \mu^{-1}$ is a BD triple on $W$ with parameter array $(\lbrace \theta_i \rbrace _{i=0}^d; \lbrace \theta'_i \rbrace _{i=0}^d; \lbrace \theta''_i \rbrace _{i=0}^d; \lbrace \rho_i \rbrace _{i=0}^d)$.  Moreover, $\lbrace \mu (V_i) \rbrace _{i=0}^d$ \\ (resp.~$\lbrace \mu(V'_i) \rbrace _{i=0}^d$) (resp.~$\lbrace \mu(V''_i) \rbrace _{i=0}^d$) is the standard ordering of the eigenspaces of $\mu A \mu^{-1}$ (resp.~$\mu A' \mu^{-1}$) (resp.~$\mu A'' \mu^{-1}$). 
\end{lemma}

\noindent
{\it Proof:}  Let $v \in V$.  Then for $0 \leq i \leq d$, the following statements are equivalent:
\begin{align*}
\mu (v) &\in \mu(V_i) \\
\iff v &\in V_i \\
\iff Av &= \theta_i v \\
\iff \mu A \mu^{-1} \mu v &= \theta_i \mu v.
\end{align*}
Thus, for $0 \leq i \leq d$, $\theta_i$ is an eigenvalue of $\mu A \mu^{-1}$ and $\mu (V_i)$ is the corresponding eigenspace.  Since $\lbrace V_i \rbrace _{i=0}^d$ is a decomposition of $V$ and $\mu: V \to W$ is an isomorphism then \\ $W = \sum_{i=0}^d \mu (V_i)$ (direct sum).  So $\mu A \mu^{-1}$ is diagonalizable.  Similarly, for $0 \leq i \leq d$, $\theta'_i$ (resp.~$\theta''_i$) is an eigenvalue of $\mu A' \mu^{-1}$ (resp.~$\mu A'' \mu^{-1}$) and $\mu (V'_i)$ (resp.~$\mu (V''_i)$) is the corresponding eigenspace.  Also $W = \sum_{i=0}^d \mu (V'_i)$ (direct sum) and $W = \sum_{i=0}^d \mu (V''_i)$ (direct sum).  So $\mu A' \mu^{-1}$ and $\mu A'' \mu^{-1}$ are diagonalizable.  By (\ref{eq:bt1}) -- (\ref{eq:bt3}) we have for $0 \leq i \leq d$,
\begin{align}
\label{muiso1}
\mu A' \mu^{-1} \mu(V_i) \subseteq \mu (V_i) + \mu (V_{i+1}), \qquad \mu A'' \mu^{-1} \mu(V_i) \subseteq \mu (V_{i-1}) + \mu (V_i), \\
\label{muiso2}
\mu A'' \mu ^{-1} \mu (V'_i) \subseteq \mu (V'_i) + \mu (V'_{i+1}),  \qquad \mu A \mu ^{-1} \mu (V'_i) \subseteq \mu (V'_{i-1}) + \mu (V'_i), \\
\label{muiso3}
\mu A \mu ^{-1} \mu (V''_i) \subseteq \mu (V''_i) + \mu (V''_{i+1}),  \qquad \mu A' \mu ^{-1} \mu (V''_i) \subseteq \mu (V''_{i-1}) + \mu (V''_i). 
\end{align}
Since $\mu :V \to W$ is a bijection then $\mu |_{V_{d-i}} : V_{d-i} \to \mu (V_{d-i})$ and $\mu ^{-1} |_{\mu (V_i)} : \mu (V_i) \to V_i$ are bijections for $0 \leq i \leq d$.  Combining this and (\ref{eq:bt4}) with Lemma \ref{thm:muj} we find that for $0 \leq i \leq d/2$, the restrictions
\begin{align*}
[\mu A' \mu ^{-1}, \mu A \mu ^{-1}]^{d-2i} |_{\mu (V_i}) : \mu (V_i) \rightarrow \mu (V_{d-i}) 
\end{align*}
are bijections.
Similarly, for $0 \leq i \leq d/2$, the restrictions
\begin{align*}
[\mu A'' \mu ^{-1}, \mu A \mu ^{-1}]^{d-2i} |_{\mu (V_{d-i})} : \mu (V_{d-i}) \rightarrow \mu (V_i), \\
[\mu A'' \mu ^{-1}, \mu A' \mu ^{-1}]^{d-2i} |_{\mu (V'_i)} : \mu (V'_i) \rightarrow \mu (V'_{d-i}), \\
[\mu A \mu ^{-1}, \mu A' \mu ^{-1}]^{d-2i} |_{\mu (V'_{d-i})} : \mu (V'_{d-i}) \rightarrow \mu (V'_i), \\
[\mu A \mu ^{-1}, \mu A'' \mu ^{-1}]^{d-2i} |_{\mu (V''_i)} : \mu (V''_i) \rightarrow \mu (V''_{d-i}), \\
[\mu A' \mu ^{-1}, \mu A'' \mu ^{-1}]^{d-2i} |_{\mu (V''_{d-i})} : \mu (V''_{d-i}) \rightarrow \mu (V''_i),
\end{align*}
are bijections.  We have now shown that $\mu A \mu^{-1}, \, \mu A' \mu^{-1}, \, \mu A'' \mu^{-1}$ is a BD triple on $W$.  By (\ref{muiso1}) (resp.~(\ref{muiso2})) (resp.~(\ref{muiso3})) we find that $\lbrace \mu (V_i) \rbrace _{i=0}^d$ (resp.~$\lbrace \mu(V'_i) \rbrace _{i=0}^d$) (resp.~$\lbrace \mu(V''_i) \rbrace _{i=0}^d$) is the standard ordering of the eigenspaces of $\mu A \mu^{-1}$ (resp.~$\mu A' \mu^{-1}$) (resp.~$\mu A'' \mu^{-1}$).  Thus, $\lbrace \theta_i \rbrace _{i=0}^d$ (resp.~$\lbrace \theta'_i \rbrace _{i=0}^d$) (resp.~$\lbrace \theta''_i \rbrace _{i=0}^d$) is the first (resp.~second) (resp.~third) eigenvalue sequence of $\mu A \mu ^{-1}, \, \mu A' \mu ^{-1}, \, \mu A'' \mu ^{-1}$.  Since $\mu : V \to W$ is a bijection then $\rho_i = {\rm dim} \, V_i = {\rm dim} \, \mu (V_i)$ ($0 \leq i \leq d$).  So $\lbrace \rho_i \rbrace _{i=0}^d$ is the shape of  $\mu A \mu ^{-1}, \, \mu A' \mu ^{-1}, \, \mu A'' \mu ^{-1}$.  We now have that $(\lbrace \theta_i \rbrace _{i=0}^d; \lbrace \theta'_i \rbrace _{i=0}^d; \lbrace \theta''_i \rbrace _{i=0}^d; \lbrace \rho_i \rbrace _{i=0}^d)$ is the parameter array of $\mu A \mu ^{-1}, \, \mu A' \mu ^{-1}, \, \mu A'' \mu ^{-1}$.  
\hfill $\Box$ \\

\section{The Proof of the Extension Theorem}

In this section we prove Theorem \ref{thm:extend}.  Throughout this section we will refer to the following assumption.

\begin{assumption}
\label{assume}
\rm
Let $A, \, A'$ denote a BD pair on $V$ of diameter $d$ and base $b$.  Let $\lbrace V_i \rbrace _{i=0}^d$ (resp.~$\lbrace V'_i \rbrace _{i=0}^d$) denote the ordering of the eigenspaces of $A$ (resp.~$A'$) satisfying (\ref{eq:bp1}) (resp. (\ref{eq:bp2})).  Let $(\lbrace \theta_i \rbrace _{i=0}^d; \lbrace \theta'_i \rbrace _{i=0}^d; \lbrace \rho_i \rbrace _{i=0}^d)$ denote the parameter array of $A, \, A'$.  
\end{assumption}

Referring to Definition \ref{def:decomp} and Assumption \ref{assume}, each of the sequences  $\lbrace V_i \rbrace _{i=0}^d$, $\lbrace V'_i \rbrace _{i=0}^d$ is a decomposition of $V$.   

\begin{definition}
\label{def:V''}
\rm
With reference to Assumption \ref{assume}, define
\begin{eqnarray*}
V''_i :=  (V_0+ \cdots +V_{d-i}) \cap (V'_0 + \cdots + V'_i) \qquad (0 \leq i \leq d).
\end{eqnarray*}
\end{definition}

\begin{lemma}
\cite[Theorem 11.2]{Funk-Neubauer13}
\label{thm:V''decomp}
With reference to Definition \ref{def:decomp} and Definition \ref{def:V''}, the sequence $\lbrace V''_i \rbrace _{i=0}^d$ is a decomposition of $V$.
\end{lemma}

\begin{lemma}
\label{thm:2cornersums}
With reference to Assumption \ref{assume}  and Definition \ref{def:V''}, the following {\rm(i)}--{\rm(ii)} hold.
\begin{enumerate}
\item[\rm (i)] $V'_0 + \cdots + V'_{d-i}= V''_0 + \cdots + V''_{d-i} \qquad (0 \leq i \leq d)$. 
\item[\rm (ii)] $ V''_i + \cdots + V''_d = V_0 + \cdots + V_{d-i} \qquad (0 \leq i \leq d)$.
\end{enumerate}
\end{lemma}

\noindent
{\it Proof:}  (i) See the proof of \cite[Corollary 11.12(i)]{Funk-Neubauer13}.  \\
(ii) Similar to (i). 
\hfill $\Box$ \\

\begin{definition}
\label{def:2A''s}
\rm
With reference to Assumption \ref{assume} and Definition \ref{def:V''}, define the following linear transformations.  
\begin{enumerate}
\item[\rm (i)] Let $\mathpzc{A''} :V \to V$ be the linear transformation such that for $0 \leq i \leq d$, $V''_i$  is an eigenspace for $\mathpzc{A''}$ with eigenvalue $\theta_i$.
\item[\rm (ii)] Let $\mathscr{A''} :V \to V$ be the linear transformation such that for $0 \leq i \leq d$, $V''_i$  is an eigenspace for $\mathscr{A''}$ with eigenvalue $\theta'_{d-i}$.
\end{enumerate}
\end{definition}

\begin{lemma}
\label{thm:2AsVV'}
With reference to Assumption \ref{assume} and Definition \ref{def:2A''s}, the following {\rm(i)}--{\rm(iv)} hold.
\begin{enumerate}
\item[\rm (i)] $(\mathpzc{A''} - \theta_i I) V'_i \subseteq V'_{i-1}$ \qquad $(0 \leq i \leq d)$.
\item[\rm (ii)] $(\mathscr{A''} - \theta'_i I) V_i \subseteq V_{i-1}$ \qquad $(0 \leq i \leq d)$.
\item[\rm (iii)] $\mathpzc{A''} \, V'_i \subseteq V'_{i-1} + V'_i$ \qquad $(0 \leq i \leq d)$.
\item[\rm (iv)] $\mathscr{A''} \, V_i \subseteq V_{i-1} + V_i$ \qquad $(0 \leq i \leq d)$.
\end{enumerate}
\end{lemma}

\noindent
{\it Proof:}  (i)  First assume that $d=0$.  In this case $V = V'_0$ and $V = V''_0$.  Thus, $(\mathpzc{A''} - \theta_0 I)V'_0 = (\mathpzc{A''} - \theta_0 I) V''_0 = 0$ and so the result holds.  Now assume that $d > 0$.  By \cite[Theorem 5.3]{Funk-Neubauer13} and \cite[Lemma 8.1]{Funk-Neubauer13} we have  
\begin{eqnarray}
\label{2As1}
\theta_{i-1} &=& b^{-1} (\theta_i - \alpha^*), \\
\label{2As2}
\theta'_{i-1} &=& b \theta'_i + \alpha, \\
\label{2As3}
\theta_{i-1} \theta'_i &=& b^{-1} (\gamma +  \theta_i \theta'_{i-1}),
\end{eqnarray}
where $\alpha, \, \alpha^*, \gamma$ denote the scalars from \cite[Theorem 5.3]{Funk-Neubauer13}.  Substituting (\ref{2As2}) into the right-hand side of (\ref{2As3}) we have
\begin{eqnarray}
\label{2As4}
\theta_{i-1} \theta'_i &=& b^{-1} \gamma + \theta_i \theta'_i + b^{-1} \alpha \theta_i.
\end{eqnarray}
By \cite[Lemma 11.4(ii)]{Funk-Neubauer13} we have $(\mathpzc{A''} - \theta_{i-1} I)(A' - \theta'_i I)V_i'' = 0$.  Substituting (\ref{2As1}) and (\ref{2As4}) into this gives 
\begin{eqnarray}
\label{2As5}
(\mathpzc{A''} A' - \theta'_i \mathpzc{A''} - b^{-1} (\theta_i - \alpha^*) A' + (b^{-1} \gamma + \theta_i \theta'_i + b^{-1} \alpha \theta_i) I) V_i'' = 0.
\end{eqnarray}
Recall that $V''_i$  is an eigenspace for $\mathpzc{A''}$ with eigenvalue $\theta_i$.  Hence (\ref{2As5}) becomes \\
($\mathpzc{A''} A' - b^{-1} A' \mathpzc{A''} + b^{-1} \alpha^* A' + b^{-1} \alpha \mathpzc{A''} + b^{-1} \gamma I) V_i'' = 0$ for $0 \leq i \leq d$.  From this and Lemma \ref{thm:V''decomp} we have 
\begin{eqnarray}
\label{2As6}
A' \mathpzc{A''} - b \mathpzc{A''} A' - \alpha \mathpzc{A''} - \alpha^* A' - \gamma I = 0.
\end{eqnarray}
Recall that $V'_i$  is an eigenspace for $A'$ with eigenvalue $\theta'_i$.  So by (\ref{2As6}) we have \\ $(A' \mathpzc{A''} - (b \theta'_i + \alpha) \mathpzc{A''} - \alpha^* \theta'_i - \gamma I) V'_i = 0$.  Substituting (\ref{2As1}) -- (\ref{2As3}) into this and simplifying gives $(A' \mathpzc{A''} - \theta'_{i-1} \mathpzc{A''} - \theta_i A' + \theta_i \theta'_{i-1} I) V'_i = 0$.  From this we have $(A' - \theta'_{i-1} I) (\mathpzc{A''} - \theta_i I) V'_i = 0$, and so $(\mathpzc{A''} - \theta_i I) V'_i \subseteq V'_{i-1}$. \\
(ii) Similar to (i). \\
(iii), (iv) Immediate from (i), (ii).
\hfill $\Box$ \\

\begin{lemma}
\label{thm:4bij}
With reference to Assumption \ref{assume}, Definition \ref{def:V''}, and Definition \ref{def:2A''s}, the following {\rm(i)}--{\rm(iv)} hold.
\begin{enumerate}
\item[\rm (i)]  For $0 \leq i \leq d/2$, the restriction $[\mathscr{A''}, A]^{d-2i} |_{V_{d-i}} : V_{d-i} \rightarrow V_i$ is a bijection.
\item[\rm (ii)]  For $0 \leq i \leq d/2$, the restriction $[\mathpzc{A''}, A']^{d-2i} |_{V'_{d-i}} : V'_{d-i} \rightarrow V'_i$ is a bijection.
\item[\rm (iii)]  For $0 \leq i \leq d/2$, the restriction $[A, \mathscr{A''}]^{d-2i} |_{V''_i} : V''_i \rightarrow V''_{d-i}$ is a bijection. 
\item[\rm (iv)]  For $0 \leq i \leq d/2$, the restriction $[A', \mathpzc{A''}]^{d-2i} |_{V''_{d-i}} : V''_{d-i} \rightarrow V''_i$ is a bijection.
\end{enumerate}
\end{lemma}

\noindent
{\it Proof:} (i) First we show that 
\begin{eqnarray}
\label{bij1}
[\mathscr{A''}, A]|_{V_i} = (\theta_i - \theta_{i-1})(\mathscr{A''} - \theta'_i I)|_{V_i} \qquad (0 \leq i \leq d).
\end{eqnarray}
Let $v \in V_i$.  By Lemma \ref{thm:2AsVV'}(ii) there exists $w \in V_{i-1}$ such that $\mathscr{A''} v = \theta'_i v + w$.  Using this we have
\begin{align*}
[\mathscr{A''}, A]v &= (\theta_i I - A) \mathscr{A''} v  \\
&= \theta'_i (\theta_i I - A)v + (\theta_i I - A)w \\
&=  (\theta_i - \theta_{i-1}) w  \\
&= (\theta_i  - \theta_{i-1}) (\mathscr{A''} - \theta'_i I) v.
\end{align*}   
We have now shown (\ref{bij1}).  Using (\ref{bij1}) and induction on $j$ we find that  
\begin{eqnarray}
\label{bij2}
[\mathscr{A''}, A]^j |_{V_i} = \Pi_{k= i+1-j}^i (\theta_k - \theta_{k-1})(\mathscr{A''} - \theta'_k I)|_{V_i} \qquad (0 \leq i \leq d) \qquad (0 \leq j \leq i+1).
\end{eqnarray}
Combining Lemma \ref{thm:raise}(ii) and Lemma \ref{thm:2AsVV'}(iv) we have that $[\mathscr{A''}, A] |_{V_i}$  maps $V_i$ into $V_{i-1}$ for $0 \leq i \leq d$.  So $[\mathscr{A''}, A]^{d-2i} |_{V_{d-i}}$  maps $V_{d-i}$ into $V_i$ for $0 \leq i \leq d/2$.  We now show that $[\mathscr{A''}, A]^{d-2i} |_{V_{d-i}}: V_{d-i} \rightarrow V_i$ is an injection for $0 \leq i \leq d/2$.  Let $ x \in {\rm ker}([\mathscr{A''}, A]^{d-2i} |_{V_{d-i}})$.  We show $x = 0$.  By construction
\begin{eqnarray}
\label{bij3}
x \in V_{d-i}.
\end{eqnarray}   
By (\ref{bij3}) and Lemma \ref{thm:2cornersums}(ii) there exist $v''_h \in V''_h$ ($i \leq h \leq d$) such that 
\begin{eqnarray}
\label{bij4}
x = \sum_{h=i}^d v''_h.
\end{eqnarray}
By construction $[\mathscr{A''}, A]^{d-2i}x =0$.  From this and (\ref{bij2}) we have $\Pi_{k= i+1}^{d-i} (\mathscr{A''} - \theta'_k I)x = 0$.  Substituting (\ref{bij4}) into this gives $\sum_{h=i}^d \Pi_{k= i+1}^{d-i} (\mathscr{A''} - \theta'_k I) v''_h= 0$.  Recall that $V''_i$  is an eigenspace for $\mathscr{A''}$ with eigenvalue $\theta'_{d-i}$.  Therefore, $\sum_{h=d-i}^d \Pi_{k= i+1}^{d-i} (\theta'_{d-h} - \theta'_k I) v''_h= 0$.  Recall, by Lemma \ref{thm:V''decomp}, that $\sum_{i=0}^d V''_i$ is a direct sum.  Thus, $v''_h = 0$ for $d-i \leq h \leq d$.  From this, (\ref{bij4}), and Lemma \ref{thm:2cornersums}(i) we find that $x \in V'_0 + \cdots + V'_{d-i-1}$.  By (\ref{bij3}) and \cite[Lemma 6.6]{Funk-Neubauer13} we have that $x \in V'_{d-i} + \cdots + V'_d$.  Combining the previous two sentences with the fact that $\sum_{i=0}^d V'_i$ is a direct sum gives $x = 0$.  We have now shown that $[\mathscr{A''}, A]^{d-2i} |_{V_{d-i}}: V_{d-i} \rightarrow V_i$ is an injection for $0 \leq i \leq d/2$.  By (\ref{eq:bp3}), ${\rm dim}(V_{d-i}) = {\rm dim}(V_i)$, and so $[\mathscr{A''}, A]^{d-2i} |_{V_{d-i}}: V_{d-i} \rightarrow V_i$ is also a surjection for $0 \leq i \leq d/2$.  Thus, $[\mathscr{A''}, A]^{d-2i} |_{V_{d-i}}: V_{d-i} \rightarrow V_i$ is a bijection for $0 \leq i \leq d/2$.  \\
(ii) Similar to (i). \\
(iii) First we show that 
\begin{eqnarray}
\label{2bij1}
[A, \mathscr{A''}]|_{V''_i} = (\theta'_{d-i} - \theta'_{d-i-1})(A - \theta_{d-i} I)|_{V''_i} \qquad (0 \leq i \leq d).
\end{eqnarray}
Let $v \in V''_i$.  By \cite[Lemma 11.4(i)]{Funk-Neubauer13} there exists $w \in V''_{i+1}$ such that $A v = \theta_{d-i} v + w$.  Using this we have
\begin{align*}
[A, \mathscr{A''}]v &= (\theta'_{d-i} I - \mathscr{A''}) A v  \\
&= \theta_{d-i} (\theta'_{d-i} I - \mathscr{A''})v + (\theta'_{d-i} I - \mathscr{A''})w \\
&=  (\theta'_{d-i} - \theta'_{d-i-1}) w  \\
&= (\theta'_{d-i}  - \theta'_{d-i-1}) (A - \theta_{d-i} I) v.
\end{align*}   
We have now shown (\ref{2bij1}).  Using (\ref{2bij1}) and induction on $j$ we find that  
\begin{eqnarray}
\label{2bij2}
[A, \mathscr{A''}]^j |_{V''_i} = \Pi_{k= i}^{j+i-1} (\theta'_{d-k} - \theta'_{d-k-1})(A - \theta_{d-k} I)|_{V''_i} \,\,\, (0 \leq i \leq d) \,\,\, (0 \leq j \leq d-i+1).
\end{eqnarray}
Combining Lemma \ref{thm:raise}(i) and \cite[Lemma 11.4(iii)]{Funk-Neubauer13} we have that $[A, \mathscr{A''}] |_{V''_i}$  maps $V''_i$ into $V''_{i+1}$ for $0 \leq i \leq d$.  So $[A, \mathscr{A''}]^{d-2i} |_{V''_i}$  maps $V''_i$ into $V''_{d-i}$ for $0 \leq i \leq d/2$.  We now show that $[A, \mathscr{A''}]^{d-2i} |_{V''_i}: V''_i \rightarrow V''_{d-i}$ is an injection for $0 \leq i \leq d/2$.  Let $ x \in {\rm ker}([A, \mathscr{A''}]^{d-2i} |_{V''_i})$.  We show $x = 0$.  By construction 
\begin{eqnarray}
\label{2bij3}
x \in V''_i.
\end{eqnarray}   
By (\ref{2bij3}) and Lemma \ref{thm:2cornersums}(ii) there exist $v_h \in V_h$ ($0 \leq h \leq d-i$) such that 
\begin{eqnarray}
\label{2bij4}
x = \sum_{h=0}^{d-i} v_h.
\end{eqnarray}
By construction $[A, \mathscr{A''}]^{d-2i}x =0$.  From this and (\ref{2bij2}) we have $\Pi_{k= i}^{d-i-1} (A - \theta_{d-k} I)x = 0$.  Substituting (\ref{2bij4}) into this gives $\sum_{h=0}^{d-i} \Pi_{k= i}^{d-i-1} (A - \theta_{d-k} I) v_h= 0$.  Recall that $V_i$  is an eigenspace for $A$ with eigenvalue $\theta_i$.  Therefore, $\sum_{h=0}^i \Pi_{k= i}^{d-i-1} (\theta_h - \theta_{d-k} I) v_h= 0$.  Recall, by Assumption \ref{assume}, that $\sum_{i=0}^d V_i$ is a direct sum.  Thus, $v_h = 0$ for $0 \leq h \leq i$.  From this, (\ref{2bij4}), and \cite[Lemma 6.6]{Funk-Neubauer13} we find that $x \in V'_{i+1} + \cdots + V'_d$.  By (\ref{2bij3}) and Lemma \ref{thm:2cornersums}(i) we have that $x \in V'_0 + \cdots + V'_i$.  Combining the previous two sentences with the fact that $\sum_{i=0}^d V'_i$ is a direct sum gives $x = 0$.  We have now shown that $[A, \mathscr{A''}]^{d-2i} |_{V''_i}: V''_i \rightarrow V''_{d-i}$ is an injection for $0 \leq i \leq d/2$.  By \cite[Corollary 11.12(ii)]{Funk-Neubauer13} and (\ref{eq:bp4}) we have ${\rm dim}(V''_i) = {\rm dim}(V''_{d-i})$, and so $[A, \mathscr{A''}]^{d-2i} |_{V''_i}: V''_i \rightarrow V''_{d-i}$ is also a surjection for $0 \leq i \leq d/2$.  Thus, $[A, \mathscr{A''}]^{d-2i} |_{V''_i}: V''_i \rightarrow V''_{d-i}$ is a bijection for $0 \leq i \leq d/2$. \\
(iv) Similar to (iii).
\hfill $\Box$ \\

We are now ready to prove the extension theorem. \\

\noindent
{\it Proof of Theorem \ref{thm:extend}:}  Adopt Assumption \ref{assume}.  Define a sequence of scalars $\lbrace \theta''_i \rbrace _{i=0}^d$ as follows.  For $d \geq 2$, let $\lbrace \theta''_i \rbrace _{i=0}^d$ be any $b$-recurrent sequence of scalars taken from $\K$.  Note that by Lemma \ref{thm:brecuraffine}(i) the sequence $\lbrace \theta''_i \rbrace _{i=0}^d$ is uniquely determined up to affine equivalence.  For $d = 1$, let $\theta''_0, \, \theta''_1$ be any two distinct scalars from $\K$.  For $d = 0$, let $\theta''_0$ be any scalar from $\K$.  With reference to Definition \ref{def:V''}, let $A'' :V \to V$ be the linear transformation such that for $0 \leq i \leq d$, $V''_i$  is an eigenspace for $A''$ with eigenvalue $\theta''_i$.  We now show that $A, \, A', \, A''$ is a BD triple on $V$.  By Definition \ref{def:bdpair}(i) and Lemma \ref{thm:V''decomp},  
\begin{eqnarray}
\label{ext1}
{\rm each \,\, of} \,\, A,\, A', \, A'' \,\, {\rm is \,\, diagonalizable}.
\end{eqnarray}
By (\ref{eq:bp1}) and (\ref{eq:bp2}),
\begin{eqnarray}
\label{ext3}
A' V_i \subseteq V_i + V_{i+1} \qquad (0 \leq i \leq d), \\
\label{ext4}
A V'_{d-i} \subseteq V'_{d-i+1} + V'_{d-i} \qquad (0 \leq i \leq d). 
\end{eqnarray}
By \cite[Lemma 11.4(iii),(iv)]{Funk-Neubauer13},
\begin{eqnarray}
\label{ext5}
A V''_i \subseteq V''_i + V''_{i+1} \qquad (0 \leq i \leq d), \\
\label{ext6}
A' V''_i \subseteq V''_{i-1} + V'_i \qquad (0 \leq i \leq d). 
\end{eqnarray}
Let $\mathpzc{A''}$ and $\mathscr{A''}$ be as in Definition \ref{def:2A''s}.  We now show that 
\begin{eqnarray}
\label{ext7}
A'' \sim \mathpzc{A''}, \\
\label{ext8}
A''  \sim \mathscr{A''}.
\end{eqnarray}
First assume that $d \geq 2$.  By \cite[Lemma 9.1]{Funk-Neubauer13} the sequences $\lbrace \theta_i \rbrace _{i=0}^d$, $\lbrace \theta'_{d-i} \rbrace _{i=0}^d$ are each $b$-recurrent, and  $\lbrace \theta''_i \rbrace _{i=0}^d$ is $b$-recurrent by construction.  From this, Lemma \ref{thm:brecuraffine}(i), and Lemma \ref{thm:mapseq} we obtain (\ref{ext7}), (\ref{ext8}).  Now assume that $d < 2$.  Then by Lemma \ref{thm:d<2} we have that $\lbrace \theta''_i \rbrace _{i=0}^d \sim \lbrace \theta_i \rbrace _{i=0}^d$ and $\lbrace \theta''_i \rbrace _{i=0}^d \sim \lbrace \theta'_{d-i} \rbrace _{i=0}^d$.  From this and Lemma \ref{thm:mapseq} we obtain (\ref{ext7}), (\ref{ext8}).    
Combining Lemma \ref{thm:2AsVV'}(iii) with (\ref{ext7}) we have
\begin{eqnarray}
\label{ext9}
A'' \, V'_{d-i} \subseteq V'_{d-i} + V'_{d-i-1} \qquad (0 \leq i \leq d).
\end{eqnarray}
Combining Lemma \ref{thm:2AsVV'}(iv) with (\ref{ext8}) we have
\begin{eqnarray}
\label{ext10}
A'' \, V_i \subseteq V_{i-1} + V_i \qquad (0 \leq i \leq d).
\end{eqnarray}
By (\ref{eq:bp3}) and (\ref{eq:bp4}),
\begin{eqnarray}
\label{ext11}
{\rm the \,\, restriction} \,\, [A', A]^{d-2i} |_{V_i} : V_i \rightarrow V_{d-i} \,\, {\rm is \,\, a \,\, bijection} \qquad (0 \leq i \leq d/2), \\
\label{ext12}
{\rm the \,\, restriction} \,\, [A, A']^{d-2i} |_{V'_i} : V'_i \rightarrow V'_{d-i} \,\, {\rm is \,\, a \,\, bijection} \qquad (0 \leq i \leq d/2).
\end{eqnarray}
Combining Lemma \ref{thm:4bij}(i),(ii) with (\ref{ext7}), (\ref{ext8}) we have
\begin{eqnarray}
\label{ext13}
{\rm the \,\, restriction} \,\, [A'', A]^{d-2i} |_{V_{d-i}} : V_{d-i} \rightarrow V_i \,\, {\rm is \,\, a \,\, bijection} \qquad (0 \leq i \leq d/2), \\
\label{ext14}
{\rm the \,\, restriction} \,\, [A'', A']^{d-2i} |_{V'_{d-i}} : V'_{d-i} \rightarrow V'_i \,\, {\rm is \,\, a \,\, bijection} \qquad (0 \leq i \leq d/2).
\end{eqnarray}
Combining Lemma \ref{thm:4bij}(iii),(iv) with (\ref{ext7}), (\ref{ext8}) we have
\begin{eqnarray}
\label{ext15}
{\rm the \,\, restriction} \,\, [A, A'']^{d-2i} |_{V''_i} : V''_i \rightarrow V''_{d-i} \,\, {\rm is \,\, a \,\, bijection} \qquad (0 \leq i \leq d/2), \\
\label{ext16}
{\rm the \,\, restriction} \,\, [A', A'']^{d-2i} |_{V''_{d-i}} : V''_{d-i} \rightarrow V''_i \,\, {\rm is \,\, a \,\, bijection} \qquad (0 \leq i \leq d/2).
\end{eqnarray}
Combining (\ref{ext1})--(\ref{ext6}) and (\ref{ext9})--(\ref{ext16}) we find that $A, \, A', \, A''$ is a BD triple on $V$.  By (\ref{ext3}) and (\ref{ext10}), $\lbrace V_i \rbrace _{i=0}^d$ is the standard ordering of the eigenspaces of $A$.  So $\lbrace \theta_i \rbrace _{i=0}^d$ is the first eigenvalue sequence of $A, \, A', \, A''$.  By (\ref{ext4}) and (\ref{ext9}), $\lbrace V'_{d-i} \rbrace _{i=0}^d$ is the standard ordering of the eigenspaces of $A'$.  So $\lbrace \theta'_{d-i} \rbrace _{i=0}^d$ is the second eigenvalue sequence of $A, \, A', \, A''$.  By Assumption \ref{assume}, $\lbrace \rho_i \rbrace _{i=0}^d$ is the shape of $A, \, A'$.  Thus, by  \cite[Definition 2.8]{Funk-Neubauer13} and Definition \ref{def:shape}, $\lbrace \rho_i \rbrace _{i=0}^d$ is the shape of $A, \, A', \, A''$.  For $d \geq 2$, \cite[Lemma 9.1] {Funk-Neubauer13} and Definition \ref{def:base} show that $A, \, A', \, A''$ has base $b$.  For $d<2$, $A, \, A'$ has base $b=1$, and so $A, \, A', \, A''$ has base $b$ by Definition \ref{def:base}.
\hfill $\Box$ \\

\section{The Proofs of the Uniqueness and Classification \\ Theorems}

In this section we prove Theorem \ref{thm:unique} and Theorem \ref{thm:class}.  First we prove the uniqueness theorem. \\

\noindent
{\it Proof of Theorem \ref{thm:unique}:}  Adopt Assumption \ref{assump}.  Let $B, \, B', \, B''$ denote a BD triple on $V$ of diameter $d$.  Let $\lbrace W_i \rbrace _{i=0}^d$ (resp.~$\lbrace W'_i \rbrace _{i=0}^d$) (resp.~$\lbrace W''_i \rbrace _{i=0}^d$) denote the standard ordering of the eigenspaces of $B$ (resp.~$B'$) (resp.~$B''$).  Let $\lbrace \sigma_i \rbrace _{i=0}^d$ (resp.~$\lbrace \sigma'_i \rbrace _{i=0}^d$) (resp.~$\lbrace \sigma''_i \rbrace _{i=0}^d$) denote the first (resp. second) (resp. third) eigenvalue sequence of $B, \, B', \, B''$.  Suppose $A \sim B$ and $A' \sim B'$.  We show that $A'' \sim B''$.  By Lemma \ref{thm:triple3pairs}(i), $A, \, A'$ is a BD pair and $\lbrace V_i \rbrace _{i=0}^d$ is an ordering of the eigenspaces of $A$ satisfying (\ref{eq:bp1}).  Since $A \sim B$ then for $0 \leq i \leq d$, $W_i$ is an eigenspace for $A$.  By (\ref{eq:bt1}), $B' W_i \subseteq W_i + W_{i+1}$ for $0 \leq i \leq d$.  Since $A' \sim B'$ then $A' W_i \subseteq W_i + W_{i+1}$ for $0 \leq i \leq d$.  Hence, $\lbrace W_i \rbrace _{i=0}^d$ is an ordering of the eigenspaces of $A$ satisfying (\ref{eq:bp1}).  Since the ordering of the eigenspaces of $A$ satisfying (\ref{eq:bp1}) is uniquely determined by $A, \, A'$ we have that 
\begin{eqnarray}
\label{uniq1}
V_i = W_i \qquad (0 \leq i \leq d).
\end{eqnarray}     
Interchanging the roles of $A$ and $A'$ and the roles of $B$ and $B'$ in the above argument gives 
\begin{eqnarray}
\label{uniq2}
V'_i = W'_i \qquad (0 \leq i \leq d).
\end{eqnarray} 
Combining (\ref{uniq1}), (\ref{uniq2}), and Lemma \ref{thm:intersect}(iii) we find that 
\begin{eqnarray}
\label{uniq3}
V''_i = W''_i \qquad (0 \leq i \leq d).
\end{eqnarray} 
We now show that 
\begin{eqnarray}
\label{uniq4}
\lbrace \theta''_i \rbrace _{i=0}^d \sim \lbrace \sigma''_i \rbrace _{i=0}^d.
\end{eqnarray}
First assume that $d \geq 2$.  Each of the sequences $\lbrace \theta_i \rbrace _{i=0}^d$, $\lbrace \theta''_i \rbrace _{i=0}^d$ is $b$-recurrent by Lemma \ref{thm:recur}.  Combining $A \sim B$ and (\ref{uniq1}) with Lemma \ref{thm:mapseq} we have that $\lbrace \theta_i \rbrace _{i=0}^d \sim \lbrace \sigma_i \rbrace _{i=0}^d$.  Hence, $\lbrace \sigma_i \rbrace _{i=0}^d$ is $b$-recurrent by Lemma \ref{thm:brecuraffine}(ii).  Therefore, $\lbrace \sigma''_i \rbrace _{i=0}^d$ is $b$-recurrent by Lemma \ref{thm:recur}.  From this and Lemma \ref{thm:brecuraffine}(i) we obtain (\ref{uniq4}).  If $d < 2$ then (\ref{uniq4}) is immediate from Lemma \ref{thm:d<2}.  Combining (\ref{uniq3}) and (\ref{uniq4}) with Lemma \ref{thm:mapseq} we find that $A'' \sim B''$.   
\hfill $\Box$ \\

The following corollary states that the linear transformation $A''$ constructed in the proof of Theorem \ref{thm:extend} is uniquely determined up to affine equivalence.

\begin{corollary}
\label{thm:extuniq}
Let $A, \, A'$ denote a BD pair on $V$.  Suppose $A''_1 : V \to V$ and $A''_2 : V \to V$ are linear transformations such that $A, \, A', \, A''_1$ and $A, \, A', \, A''_2$ are both BD triples on $V$.  Then $A''_1 \sim A''_2$.  
\end{corollary}

\noindent
{\it Proof:}  Immediate from Theorem \ref{thm:unique}.
\hfill $\Box$ \\

Lemma \ref{thm:isompa} will be used in the proof of Theorem \ref{thm:class}.  So we now prove Lemma \ref{thm:isompa}, namely that two BD triples are isomorphic exactly when their parameter arrays are equal. \\

\noindent
{\it Proof of Lemma \ref{thm:isompa}:}  Let $A, \, A', \, A''$ and $B, \, B', \, B''$ denote BD triples over $\K$.  Let $V$ (resp.~$W$) denote the vector space underlying $A, \, A', \, A''$ (resp.~$B, \, B', \, B''$).  Let
\begin{eqnarray}
\label{isompa1}
(\{ \theta_i \}_{i=0}^d; \{  \theta'_i \}_{i=0}^d; \{  \theta''_i \}_{i=0}^d; \{ \rho_i \}_{i=0}^d)
\end{eqnarray}
denote the parameter array of $A, \, A', \, A''$.  \\
$(\Longrightarrow)$:  Suppose that $A, \, A', \, A''$ and $B, \, B', \, B''$ are isomorphic.  We show that (\ref{isompa1}) is also the parameter array of $B, \, B', \, B''$.  Let $\mu : V \rightarrow W$ denote an isomorphism of BD triples from $A, \, A', \, A''$ to $B, \, B', \, B''$.  By Lemma \ref{thm:triple3pairs}(i), $A, \, A'$ is a BD pair on $V$ with parameter array $(\{ \theta_i \}_{i=0}^d; \{  \theta'_{d-i} \}_{i=0}^d; \{ \rho_i \}_{i=0}^d)$, and $B, \, B'$ is a BD pair on $W$.  Comparing \cite[Definition 2.11]{Funk-Neubauer13} and Definition \ref{def:isom} we find that $\mu : V \rightarrow W$ is an isomorphism of BD pairs from $A, \, A'$ to $B, \, B'$.  So by \cite[Lemma 2.12]{Funk-Neubauer13}, $(\{ \theta_i \}_{i=0}^d; \{  \theta'_{d-i} \}_{i=0}^d; \{ \rho_i \}_{i=0}^d)$ is the parameter array of $B, \, B'$.  Hence,  $\{ \rho_i \}_{i=0}^d$ is the shape of $B, \, B'$, and so $\{ \rho_i \}_{i=0}^d$ is the shape of $B, \, B', \, B''$.  Moreover, $\{ \theta_i \}_{i=0}^d$ is the eigenvalue sequence of $B, \, B'$, and so $\{ \theta_i \}_{i=0}^d$ is the first eigenvalue sequence of $B, \, B', \, B''$.  Similarly, $\{ \theta'_i \}_{i=0}^d$ (resp.~$\{ \theta''_i \}_{i=0}^d$) is the second (resp. third) eigenvalue sequence of $B, \, B', \, B''$.  Combining the previous three sentences we have that (\ref{isompa1}) is the parameter array of $B, \, B', \, B''$.  Therefore, the parameter array of $A, \, A', \, A''$ equals the parameter array of $B, \, B', \, B''$. \\    
$(\Longleftarrow)$:  Suppose that the parameter array of $A, \, A', \, A''$ equals the parameter array of $B, \, B', \, B''$.  That is, suppose  
\begin{eqnarray}
\label{isompa1.5}
{\rm (\ref{isompa1}) \,\, is \,\, also \,\, the \,\, parameter \,\, array \,\, of \,\,} B, \, B', \, B''.  
\end{eqnarray}
We show that $A, \, A', \, A''$ and $B, \, B', \, B''$ are isomorphic.  By Lemma \ref{thm:triple3pairs}(i), the BD pairs $A, \, A'$  and $B, \, B'$ have the same parameter array.  So by \cite[Lemma 2.12]{Funk-Neubauer13}, $A, \, A'$ and $B, \, B'$ are isomorphic.  So there exists a vector space isomorphism $\mu : V \to W$ such that 
\begin{eqnarray}
\label{isompa2}
\mu A = B \mu, \qquad \mu A' = B' \mu.
\end{eqnarray} 
We now show that 
\begin{eqnarray}
\label{isompa2.5}
\mu A'' = B'' \mu.
\end{eqnarray}  
By Lemma \ref{thm:muiso}, $\mu A \mu^{-1}, \, \mu A' \mu^{-1}, \, \mu A'' \mu^{-1}$ is a BD triple on $W$.  So by Corollary \ref{thm:extuniq} and (\ref{isompa2}), we have $\mu A'' \mu ^{-1} \sim B''$.  Thus, there exist $r, s$ in $\K$ such that  
\begin{eqnarray}
\label{isompa3}
\mu A'' \mu ^{-1} = rB'' + sI.
\end{eqnarray}
Let $\{ V''_i \}_{i=0}^d$ denote the standard ordering of the eigenspaces of $A''$.  By Lemma \ref{thm:muiso}, $\lbrace \mu(V''_i) \rbrace _{i=0}^d$ is the standard ordering of the eigenspaces of $\mu A'' \mu^{-1}$.  From this, (\ref{isompa2}), and (\ref{isompa3}) we find that $\lbrace \mu(V''_i) \rbrace _{i=0}^d$ is also the standard ordering of the eigenspaces of $B''$.  By Lemma \ref{thm:muiso} and (\ref{isompa1.5}), $\{  \theta''_i \}_{i=0}^d$ is the third eigenvalue sequence of both $\mu A \mu^{-1}, \, \mu A' \mu^{-1}, \, \mu A'' \mu^{-1}$ and $B, \, B', \, B''$.  We have now shown that for $0 \leq i \leq d$, $\mu (V''_i)$ is an eigenspace for both $\mu A'' \mu ^{-1}$ and $B''$ with eigenvalue $\theta''_i$.  Combining this and (\ref{isompa3}) with Lemma \ref{thm:mapseq} we have 
\begin{eqnarray}
\label{isompa4}
\theta''_i = r \theta''_i + s \qquad (0 \leq i \leq d).  
\end{eqnarray}
Now assume that $d \geq 1$.  Subtracting (\ref{isompa4}) with $i=1$ from (\ref{isompa4}) with $i=0$ gives $\theta''_0 - \theta''_1 = r(\theta''_0 - \theta''_1)$.  Since $\theta''_0 \neq \theta''_1$ we find that $r=1$.  Substituting this into (\ref{isompa4}) gives $s=0$.  Substituting $r=1$ and $s=0$ into (\ref{isompa3}) gives (\ref{isompa2.5}).  Next assume that $d < 1$.  In this case $\mu A'' \mu^{-1}$ and $B''$ are both equal to $\theta''_0 I$, and so (\ref{isompa2.5}) holds.  Combining (\ref{isompa2}) and (\ref{isompa2.5}) we find that $\mu$ is an isomorphism of bidiagonal triples from $A, \, A', \, A''$ to $B, \, B', \, B''$.  So $A, \, A, \, A''$ and $B, \, B, \, B''$ are isomorphic.
\hfill $\Box$ \\

\noindent
We are now ready to prove the classification theorem.  \\

\noindent
{\it Proof of Theorem \ref{thm:class}:} \\
$(\Longrightarrow)$:  Adopt Assumption \ref{assump}.  We show that the scalars in (\ref{parameterarray}) satisfy Theorem \ref{thm:class}(i)--(v).  \\ 
First assume that $d < 1$.  Theorem \ref{thm:class}(i),(ii),(v) hold since they are vacuously true.  Theorem \ref{thm:class}(iv) holds since it is trivially true.  Recall $V_0$ is an eigenspace of $A$ and $\rho_0 = \dim(V_0)$ by Definition \ref{def:shape}.  Thus, $\rho_0$ is a positive integer and Theorem \ref{thm:class}(iii) holds.  Now assume that $d \geq 1$.  By definition of first (resp.~second) (resp.~third) eigenvalue sequence $\{ \theta_i \}_{i=0}^d$ (resp.~$\{ \theta'_i \}_{i=0}^d$) (resp.~$\{ \theta''_i \}_{i=0}^d$) is a list of the {\it distinct} eigenvalues of $A$ (resp.~$A'$)  (resp.~$A''$).  Thus, Theorem \ref{thm:class}(i) holds.  For $d=1$, Theorem \ref{thm:class}(ii) holds since it is vacuously true.  For $d > 1$, Theorem \ref{thm:class}(ii) holds by Lemma \ref{thm:recur}.  Recall for $0 \leq i \leq d$, $V_i$ is an eigenspace of $A$ and $\rho_i = \dim(V_i)$ by Definition \ref{def:shape}.  Thus, $\rho_i$ is a positive integer for $0 \leq i \leq d$ and Theorem \ref{thm:class}(iii) holds.  By Definition \ref{def:shape}, $\rho_i = \dim(V_i) = \dim(V_{d-i}) = \rho_{d-i}$ for $0 \leq i \leq d$.  So Theorem \ref{thm:class}(iv) holds.  By (\ref{eq:bt4}), the restriction $[A',A] |_{V_i} : V_i \rightarrow V_{i+1}$ is an injection for $0 \leq i < d/2$.  Thus, by Definition \ref{def:shape}, $\rho_i = \dim(V_i) \leq \dim(V_{i+1}) = \rho_{i+1}$ for $0 \leq i < d/2$.  So Theorem \ref{thm:class}(v) holds. \\ 
$(\Longleftarrow)$:  Let (\ref{parameterarray}) denote a sequence of scalars taken from $\K$ which satisfy Theorem \ref{thm:class}(i)--(v).  We construct a BD triple $A, \, A', \, A''$ over $\K$ which has parameter array (\ref{parameterarray}).  Observe that $(\lbrace \theta_i \rbrace _{i=0}^d; \lbrace \theta'_{d-i} \rbrace _{i=0}^d; \lbrace \rho_i \rbrace _{i=0}^d)$ satisfy \cite[Theorem 5.1(i)--(v)]{Funk-Neubauer13}.  Therefore, by \cite[Theorem 5.1]{Funk-Neubauer13}, there exists a BD pair $A, \, A'$ over $\K$ with parameter array $(\lbrace \theta_i \rbrace _{i=0}^d; \lbrace \theta'_{d-i} \rbrace _{i=0}^d; \lbrace \rho_i \rbrace _{i=0}^d)$.  Let $V$ denote the vector space underlying $A, \, A'$, and let $b$ denote the base of $A, \, A'$.  By Theorem \ref{thm:extend} there exists a linear transformation $\overline{A''}: V \rightarrow V$ such that $A, \, A', \, \overline{A''}$ is a BD triple on $V$ of base $b$.  Let $\lbrace V''_i \rbrace _{i=0}^d$ denote the standard ordering of eigenspaces of $\overline{A''}$.  Let $\lbrace \overline{\theta''_i} \rbrace _{i=0}^d$ denote the third eigenvalue sequence of $A, \, A', \, \overline{A''}$.  By Theorem \ref{thm:extend} 
\begin{eqnarray}
\label{class1}
A, \, A', \, \overline{A''} {\rm \,\, has \,\, parameter \,\, array} \,\, (\lbrace \theta_i \rbrace _{i=0}^d; \lbrace \theta'_i \rbrace _{i=0}^d; \lbrace \overline{\theta''_i} \rbrace _{i=0}^d; \lbrace \rho_i \rbrace _{i=0}^d).
\end{eqnarray}   
We now show that 
\begin{eqnarray}
\label{class2}
\lbrace \theta''_i \rbrace _{i=0}^d \sim \lbrace \overline{\theta''_i} \rbrace _{i=0}^d.
\end{eqnarray}
First assume that $d \geq 2$.  By Lemma \ref{thm:recur} and Theorem \ref{thm:class}(ii), each of the sequences $\lbrace \theta''_i \rbrace _{i=0}^d$, $\lbrace \overline{\theta''_i} \rbrace _{i=0}^d$ is $b$-recurrent.  From this and Lemma \ref{thm:brecuraffine}(i) we obtain (\ref{class2}).  If $d < 2$ then (\ref{class2}) is immediate from Lemma \ref{thm:d<2}.  Let $A'': V \to V$ be the linear transformation such that for $0 \leq i \leq d$, $V''_i$  is an eigenspace for $A''$ with eigenvalue $\theta''_i$.  Combining (\ref{class2}) and Lemma \ref{thm:mapseq} we find that 
$A'' \sim \overline{A''}$.  Thus, by Lemma \ref{thm:affineshift} and (\ref{class1}) we have that $A, \, A', \, A''$ is a BD triple over $\K$ with parameter array (\ref{parameterarray}).  The uniqueness claim in Theorem \ref{thm:class} follows immediately from Lemma \ref{thm:isompa}.
\hfill $\Box$ \\

\section{The Proofs of the Remaining Theorems}

In this section we prove Theorem \ref{thm:triplerelations}, Theorem \ref{thm:reduce}, Theorem \ref{thm:uqsl2bidiag}, and Theorem \ref{thm:reducedtomod}.  The following lemma will be used in the proof of Theorem \ref{thm:triplerelations}. 

\begin{lemma}
\label{thm:6recurs}
With reference to Assumption \ref{assump} suppose that $d \geq 1$.  Then there exists a sequence of scalars $\alpha, \, \alpha', \, \alpha'', \, \gamma_1, \, \gamma_2, \, \gamma_3$ in $\K$ such that the following {\rm (i)--(vi)} hold.
\begin{enumerate}
\item[\rm (i)] $\theta_{i+1} - b \theta_i = \alpha \qquad (0 \leq i \leq d-1)$,
\item[\rm (ii)] $\theta'_{i+1} - b \theta'_i = \alpha' \qquad (0 \leq i \leq d-1)$,
\item[\rm (iii)] $\theta''_{i+1} - b \theta''_i = \alpha'' \qquad (0 \leq i \leq d-1)$,
\item[\rm (iv)] $b \theta_i \theta'_{d-i-1} - \theta_{i+1} \theta'_{d-i} = \gamma_1 \qquad (0 \leq i \leq d-1)$, 
\item[\rm (v)] $b \theta'_i \theta''_{d-i-1} - \theta'_{i+1} \theta''_{d-i} = \gamma_2 \qquad (0 \leq i \leq d-1)$, 
\item[\rm (vi)] $b \theta''_i \theta_{d-i-1} - \theta''_{i+1} \theta_{d-i} = \gamma_3 \qquad (0 \leq i \leq d-1)$.
\end{enumerate}
\end{lemma}

\noindent
{\it Proof:}  For $d \geq 2$, Lemma \ref{thm:recur} and Definition \ref{def:base} give that 
\begin{eqnarray}
\label{6recurs1}
\frac{\theta_{i+1} - \theta_i}{\theta_i - \theta_{i-1}} = b = \frac{\theta'_{i+1} - \theta'_i}{\theta'_i - \theta'_{i-1}}, \qquad (1 \leq i \leq d-1).
\end{eqnarray}
(i) Define $\alpha := \theta_1 - b \theta_0$.  Then for $d=1$ the result is immediate.  For $d \geq 2$, (\ref{6recurs1}) becomes $\theta_{i+1} - b \theta_i = \theta_i - b \theta_{i-1}$ $(1 \leq i \leq d-1)$, and the result follows.  \\
(ii), (iii) Similar to (i).  \\
(iv) For $d=1$, define $\gamma_1 := b \theta_0 \theta'_0 - \theta_1 \theta'_1$ and the result is immediate.  Now assume that $d \geq 2$.  First suppose that $b=1$.  By solving the recurrences in (\ref{6recurs1}) we find that there exist $c_1, c_2, c_3, c_4$ in $\K$ such that $\theta_i = c_1 + c_2 i$ and $\theta'_i = c_3 + c_4 i$ ($0 \leq i \leq d$).  Define $\gamma_1 := - (c_1 c_4 + c_2 c_3 + c_2 c_4 d)$.  Combining the previous two sentences we obtain the result.  Now suppose that $b \neq 1$.  By solving the recurrences in (\ref{6recurs1}) we find that there exist $c_1, c_2, c_3, c_4$ in $\K$ such that $\theta_i = c_1 + c_2 b^i$ and $\theta'_i = c_3 + c_4 b^i$ ($0 \leq i \leq d$).  Define $\gamma_1 : = (b-1)(c_1 c_3 - c_2 c_4 b^d)$.  Combining the previous two sentences we obtain the result. \\
(v), (vi) Similar to (iv). 
\hfill $\Box$ \\

\noindent
We now prove that the fundamental bidiagonal relations hold. \\

\noindent
{\it Proof of Theorem \ref{thm:triplerelations}:}  Adopt Assumption \ref{assump}.  First assume that $d \geq 1$.  Let $\alpha, \, \alpha', \, \alpha'', \gamma_1, \, \gamma_2, \\ \gamma_3$ denote the scalars from Lemma \ref{thm:6recurs}.  By Lemma \ref{thm:triple3pairs}(i), $A, \, A'$ is a BD pair with parameter array $(\lbrace \theta_i \rbrace _{i=0}^d; \lbrace \theta'_{d-i} \rbrace _{i=0}^d; \lbrace \rho_i \rbrace _{i=0}^d)$.  Combining Lemma \ref{thm:6recurs}(i),(ii),(iv) with \cite[Lemma 8.1]{Funk-Neubauer13} we obtain (\ref{eq:triple1}).  By Lemma \ref{thm:triple3pairs}(ii), $A', \, A''$ is a BD pair with parameter array \\ $(\lbrace \theta'_i \rbrace _{i=0}^d; \lbrace \theta''_{d-i} \rbrace _{i=0}^d; \lbrace \rho_i \rbrace _{i=0}^d)$.  Combining Lemma \ref{thm:6recurs}(ii),(iii),(v) with \cite[Lemma 8.1]{Funk-Neubauer13} we obtain (\ref{eq:triple2}).  By Lemma \ref{thm:triple3pairs}(iii), $A'', \, A$ is a BD pair with parameter array $(\lbrace \theta''_i \rbrace _{i=0}^d; \lbrace \theta_{d-i} \rbrace _{i=0}^d; \lbrace \rho_i \rbrace _{i=0}^d)$.  Combining Lemma \ref{thm:6recurs}(i),(iii),(vi) with \cite[Lemma 8.1]{Funk-Neubauer13} we obtain (\ref{eq:triple3}).  Next we demonstrate the uniqueness claim for $d \geq 2$.  Recall that each of the sequences $\lbrace \theta_i \rbrace _{i=0}^d$, $\lbrace \theta'_i \rbrace _{i=0}^d$, $\lbrace \theta''_i \rbrace _{i=0}^d$ is uniquely determined by $A, \, A', \, A''$.  From this and Definition \ref{def:base}, $b$ is uniquely determined by $A, \, A', \, A''$.  Combining the previous two sentences with Lemma \ref{thm:6recurs}, each of $\alpha, \, \alpha', \, \alpha'', \, \gamma_1, \, \gamma_2, \, \gamma_3$ is uniquely determined by $A, \, A', \, A''$.  Now assume that $d = 0$.  Recall by Definition \ref{def:base} that $b=1$.  Define $\alpha := 1$, $\alpha' := 1$, $\alpha'' := 1$, $\gamma_1 := -(\theta_0 + \theta'_0)$, $\gamma_2 :=  -(\theta'_0 + \theta''_0)$, and $\gamma_3 := -(\theta''_0 + \theta_0)$.  Observe that $A = \theta_0 I$, $A' = \theta'_0 I$, and $A'' = \theta''_0 I$.  Combining the previous three sentences we obtain (\ref{eq:triple1})--(\ref{eq:triple3}).      
\hfill $\Box$ \\

\begin{corollary}
\label{thm:reducedrels}
With reference to Assumption \ref{assump} suppose that $A, \, A', \, A''$ is reduced.  Then the following {\rm (i),(ii)} hold.   
\begin{enumerate}
\item[\rm (i)] Suppose $b = 1$.  Then the relations (\ref{eq:triple1})--(\ref{eq:triple3}) take the form
\begin{eqnarray*}
\label{RR1}
A A' - A' A - 2A - 2 A' = 0, \\
\label{RR2}
A' A'' - A'' A' - 2A' - 2 A'' = 0, \\
\label{RR3}
A'' A - A A'' - 2A'' - 2 A = 0.
\end{eqnarray*}
\item[\rm (ii)] Suppose $b \neq 1$.  Then the relations (\ref{eq:triple1})--(\ref{eq:triple3}) take the form
\begin{eqnarray*}
\label{RR4}
q A A' - q^{-1} A' A - (q - q^{-1}) I =  0, \\
\label{RR5}
q A' A'' - q^{-1} A'' A' - (q - q^{-1}) I =  0, \\
\label{RR6}
q A'' A - q^{-1} A A'' - (q - q^{-1}) I =  0.
\end{eqnarray*}
\end{enumerate}
\end{corollary}

\noindent
{\it Proof:}  (i) First assume that $d \geq 1$.  Combining Definition \ref{def:reduced}(i) with Lemma \ref{thm:6recurs} we find that $\alpha = 2, \alpha' = 2, \alpha'' = 2, \gamma_1 = 0, \gamma_2 = 0, \gamma_3 = 0$.  Substituting these into (\ref{eq:triple1})--(\ref{eq:triple3}) we obtain the result.  Now assume that $d = 0$.  Observe that $A = \theta_0 I$, $A' = \theta'_0 I$, and $A'' = \theta''_0 I$.  By Definition \ref{def:reduced}(i), $\theta_0 =0$, $\theta'_0 = 0$, and $\theta''_0 = 0$.  Combining the previous two sentences we obtain the result.   \\
(ii)  First assume that $d \geq 1$.  Combining Definition \ref{def:reduced}(ii) with Lemma \ref{thm:6recurs} we find that $\alpha = 0, \alpha' = 0, \alpha'' = 0, \gamma_1 = q^{-1}(q-q^{-1}), \gamma_2 = q^{-1}(q-q^{-1}), \gamma_3 = q^{-1}(q-q^{-1})$.  Substituting these into (\ref{eq:triple1})--(\ref{eq:triple3}) we obtain the result.  Now assume that $d = 0$.  Observe that $A = \theta_0 I$, $A' = \theta'_0 I$, and $A'' = \theta''_0 I$.  By Definition \ref{def:reduced}(ii), $\theta_0 =1$, $\theta'_0 = 1$, and $\theta''_0 = 1$.  Combining the previous two sentences we obtain the result.          
\hfill $\Box$ \\

\noindent
We now prove the reducibility theorem. \\

\noindent
{\it Proof of Theorem \ref{thm:reduce}:}  Adopt Assumption \ref{assump}.  Suppose that $b = 1$.  First we show that 
\begin{eqnarray}
\label{reduce1}
{\rm \,\, each \,\, of \,\,} \lbrace \theta_i \rbrace _{i=0}^d, \lbrace \theta'_i \rbrace _{i=0}^d, \lbrace \theta''_i \rbrace _{i=0}^d {\rm \,\, is \,\, affine \,\, equivalent \,\, to \,\,} \lbrace 2i-d \rbrace _{i=0}^d.
\end{eqnarray}
Assume that $d \geq 2$.  By Lemma \ref{thm:recur} and Definition \ref{def:base}, each of the sequences $\lbrace \theta_i \rbrace _{i=0}^d$, $\lbrace \theta'_i \rbrace _{i=0}^d$, $\lbrace \theta''_i \rbrace _{i=0}^d$ is $b$-recurrent.  Observe that the sequence $\lbrace 2i-d \rbrace _{i=0}^d$ is also $b$-recurrent.  So by Lemma \ref{thm:brecuraffine}(i) we obtain (\ref{reduce1}).  If $d < 2$ then (\ref{reduce1}) is immediate from Lemma \ref{thm:d<2}.  Let $B :V \to V$ (resp.~$B' :V \to V$) (resp.~$B'' :V \to V$) be the linear transformation such that for $0 \leq i \leq d$, $V_i$ (resp.~$V'_i$) (resp.~$V''_i$) is an eigenspace for $B$ (resp.~$B'$) (resp.~$B''$) with eigenvalue $2i-d$.  From (\ref{reduce1}) and Lemma \ref{thm:mapseq} we have $A \sim B$, $A' \sim B'$, and $A'' \sim B''$.  Combining this with Lemma \ref{thm:affineshift} we find that $B, \, B', \, B''$ is a BD triple.  Also $A, \, A', \, A''$ is affine equivalent to $B, \, B', \, B''$.  Lemma \ref{thm:affineshift} also shows that each of the first, second, and third eigenvalue sequences of $B, \, B', \, B''$ is $\lbrace 2i-d \rbrace _{i=0}^d$.  Hence, $B, \, B', \, B''$ is reduced.  So $A, \, A', \, A''$ is affine equivalent to a reduced BD triple.  Repeating the above argument with $\lbrace 2i-d \rbrace _{i=0}^d$ replaced by $\lbrace q^{d-2i} \rbrace _{i=0}^d$ proves the result for $b \neq 1$.    
\hfill $\Box$ \\

We now prove the correspondence theorems. \\

\noindent
{\it Proof of Theorem \ref{thm:uqsl2bidiag}:}  Let $V$ denote a segregated $\SL$-module.  Let $X, \, Y, \, Z$ denote an equitable triple in $\SL$.  We show that the action of $X, \, Y, \, Z$ on $V$ is a reduced BD triple of base $1$.  By \cite[Theorem 5.10]{Funk-Neubauer13} the action of $Y, \, Z$ on $V$ is a reduced BD pair of base $1$.  By Theorem \ref{thm:extend} there exists a linear transformation $\overline{\mathcal{X}}: V \rightarrow V$ such that $Y, \, Z, \, \overline{\mathcal{X}}$ is a BD triple on $V$ of base $1$.  So $\overline{\mathcal{X}}, \, Y, \, Z$ is a BD triple on $V$ of base $1$.  Since the action of $Y, \, Z$ on $V$ is a {\it reduced} BD pair of base $1$ then the second and third eigenvalue sequences of $\overline{\mathcal{X}}, \, Y, \, Z$ are each $\lbrace 2i-d \rbrace _{i=0}^d$.  Let $\lbrace \theta_i \rbrace _{i=0}^d$ denote the first eigenvalue sequence of $\overline{\mathcal{X}}, \, Y, \, Z$.  We now show that  
\begin{eqnarray}
\label{uqsl2bidiag1}
\lbrace \theta_i \rbrace _{i=0}^d \sim \lbrace 2i-d \rbrace _{i=0}^d.  
\end{eqnarray}
First assume that $d \geq 2$.  By Lemma \ref{thm:recur} and Definition \ref{def:base} the sequence $\lbrace \theta_i \rbrace _{i=0}^d$ is $1$-recurrent.  Observe that the sequence $\lbrace 2i-d \rbrace _{i=0}^d$ is also $1$-recurrent.  So by Lemma \ref{thm:brecuraffine}(i) we obtain (\ref{uqsl2bidiag1}).  If $d < 2$ then (\ref{uqsl2bidiag1}) is immediate from Lemma \ref{thm:d<2}.  Let  $\lbrace V_i \rbrace _{i=0}^d$ denote the standard ordering of the eigenspaces of $\overline{\mathcal{X}}$.  Let $\mathcal{X}: V \to V$ be the linear transformation such that for $0 \leq i \leq d$, $V_i$ is an eigenspace for $\mathcal{X}$ with eigenvalue $2i-d$.  Combining (\ref{uqsl2bidiag1}) with Lemma \ref{thm:mapseq} we have $\overline{\mathcal{X}} \sim \mathcal{X}$.  So, by Lemma \ref{thm:affineshift}, $\mathcal{X}, \, Y, \, Z$ is a BD triple on $V$ of base $1$ whose first, second, and third eigenvalue sequences are each $\lbrace 2i-d \rbrace _{i=0}^d$.  So $\mathcal{X}, \, Y, \, Z$ is reduced.  We now show that the action of $\mathcal{X}$ on $V$ is equal to the action of $X$ on $V$.  By construction the action of  $X, \, Y, \, Z$ defines an $\SL$-module structure on $V$.  Comparing Theorem \ref{thm:ueq} with Corollary \ref{thm:reducedrels} we find that the action of $\mathcal{X}, \, Y, \, Z$ defines another $\SL$-module structure on $V$.  So by \cite[Lemma 12.5]{Funk-Neubauer13} the action of $\mathcal{X}$ on $V$ is equal to the action of $X$ on $V$.  We have now shown that the action of $X, \, Y, \, Z$ on $V$ is a reduced BD triple of base $1$.  Repeating the above argument with $\lbrace 2i-d \rbrace _{i=0}^d$ replaced by $\lbrace q^{d-2i} \rbrace _{i=0}^d$ proves the result in the $\uq$ case.   
\hfill $\Box$ \\ 

\noindent
{\it Proof of Theorem \ref{thm:reducedtomod}:}  Let $A, \, A', \, A''$ denote a reduced BD triple on $V$  of base $b$.   \\
(i) Suppose that $b = 1$.  Let $X, \, Y, \, Z$ denote an equitable triple in $\SL$.  Comparing Theorem \ref{thm:ueq} with Corollary \ref{thm:reducedrels}(i) we find that there exists an $\SL$-module structure on $V$ such that $(X - A)V=0$, $(Y - A')V=0$, and $(Z - A'')V=0$.  We now check that this $\SL$-module structure on $V$ is segregated.  Using the $\SL$ automorphism from Theorem \ref{thm:ueq} we find that $(h - Z)V = 0$.  So  $(h - A'')V = 0$.  Since $A, \, A', \, A''$ is reduced then the third eigenvalue sequence of $A, \, A', \, A''$ is $\lbrace 2i-d \rbrace _{i=0}^d$.  The previous two sentences show that the action of $h$ on $V$ has eigenvalues $\lbrace 2i-d \rbrace _{i=0}^d$.  So if $d$ is even (resp.~odd) then $V=V_{\rm{even}}$ (resp.~$V=V_{\rm{odd}}$).  Thus, the $\SL$-module structure on $V$ is segregated.   \\
(ii) Suppose that $b \neq 1$.  Let $x, \, y, \, z$ denote an equitable triple in $\uq$.  Comparing Theorem \ref{thm:uqeq} with Corollary \ref{thm:reducedrels}(ii) we find that there exists a $\uq$-module structure on $V$ such that $(x - A)V=0$, $(y - A')V=0$, and $(z - A'')V=0$.  We now check that this $\uq$-module structure on $V$ is segregated.  Using the $\uq$ automorphism from Theorem \ref{thm:uqeq} we find that $(k - x)V = 0$.  So  $(k - A)V = 0$.  Since $A, \, A', \, A''$ is reduced then the first eigenvalue sequence of $A, \, A', \, A''$ is $\lbrace q^{d-2i} \rbrace _{i=0}^d$.  The previous two sentences show that the action of $k$ on $V$ has eigenvalues $\lbrace q^{d-2i} \rbrace _{i=0}^d$.  So if $d$ is even (resp.~odd) then $V=V_{\rm{even}}^1$ (resp.~$V=V_{\rm{odd}}^1$).  Thus, the $\uq$-module structure on $V$ is segregated.     
\hfill $\Box$ \\

\section{Acknowledgment}
The author would like to thank Professor Paul Terwilliger for his advice and many helpful suggestions.  His comments greatly improved this paper.   
 
\bibliographystyle{amsplain}
\bibliography{FunkNeubauer}

\providecommand{\bysame}{\leavevmode\hbox to3em{\hrulefill}\thinspace}
\providecommand{\MR}{\relax\ifhmode\unskip\space\fi MR }
\providecommand{\MRhref}[2]{%
  \href{http://www.ams.org/mathscinet-getitem?mr=#1}{#2}
}
\providecommand{\href}[2]{#2}
\begin{thebibliography}{10}

\bibitem{Al-Najjar10}
H.~Al-Najjar and B.~Curtin, \emph{Leonard pairs from the equitable basis of
  $\mathfrak{sl}_2$}, Electron. J. Linear Algebra \textbf{20} (2010), 490--505.

\bibitem{Benkart04}
G.~Benkart and P.~Terwilliger, \emph{Irreducible modules for the quantum affine
  algebra ${U}_q(\widehat{\mathfrak{sl}}_2)$ and its {B}orel subalgebra}, J.
  Algebra \textbf{282} (2004), 172--194.

\bibitem{Benkart07}
\bysame, \emph{The universal central extension of the three-point
  $\mathfrak{sl}_2$ loop algebra}, Proc. Amer. Math. Soc. \textbf{135} (2007),
  no.~6, 1659--1668.

\bibitem{Benkart10}
\bysame, \emph{The equitable basis for $\mathfrak{sl}_2$}, Mathematische
  Zeitschrift (2010), published online at www.springerlink.com.

\bibitem{Elduque07}
A.~Elduque, \emph{The ${S}_4$-action on the tetrahedron algebra}, Proc. Roy.
  Soc. Edinburgh Sect. A \textbf{137} (2007), no.~6, 1227--1248.

\bibitem{Funk-Neubauer07}
D.~Funk-Neubauer, \emph{Raising/lowering maps and modules for the quantum
  affine algebra ${U}_q(\widehat {\mathfrak{sl}}_2)$}, Comm Algebra \textbf{35}
  (2007), no.~7, 2140--2159.

\bibitem{Funk-Neubauer09}
\bysame, \emph{Tridiagonal pairs and the $q$-tetrahedron algebra}, Linear
  Algebra Appl. \textbf{431} (2009), no.~5--7, 903--925.

\bibitem{Funk-Neubauer13}
\bysame, \emph{Bidiagonal pairs, the {L}ie algebra $\mathfrak{sl}_2$, and the
  quantum group ${U}_q(\mathfrak{sl}_2)$}, J. Algebra Appl. \textbf{12} (2013),
  no.~5, 1250207, pp. 46.

\bibitem{Hartwig07}
B.~Hartwig, \emph{The tetrahedron algebra and its finite dimensional
  irreducible modules}, Linear Algebra Appl. \textbf{422} (2007), no.~1,
  219--235.

\bibitem{HarTer07}
B.~Hartwig and P.~Terwilliger, \emph{The tetrahedron algebra, the {O}nsager
  algebra, and the $\mathfrak{sl}_2$ loop algebra}, J. Algebra \textbf{308}
  (2007), no.~2, 840--863.

\bibitem{Humph72}
J.~E. Humphreys, \emph{Introduction to {L}ie algebras and representation
  theory}, Springer-Verlag, New York, NY, 1972.

\bibitem{ItoTanTer01}
T.~Ito, K.~Tanabe, and P.~Terwilliger, \emph{Some algebra related to {P}- and
  {Q}-polynomial association schemes}, DIMACS Ser. Discrete Math. Theoret.
  Comput. Sci. (Providence, RI), vol.~56, American Mathematical Society, 2001,
  pp.~167--192.

\bibitem{ItoTer071}
T.~Ito and P.~Terwilliger, \emph{$q$-{I}nverting pairs of linear
  transformations and the $q$-tetrahedron algebra}, Linear Algebra Appl.
  \textbf{426} (2007), no.~2--3, 516--532.

\bibitem{ItoTer072}
\bysame, \emph{The $q$-tetrahedron algebra and its finite dimensional
  irreducible modules}, Comm. Algebra \textbf{35} (2007), no.~11, 3415--3439.

\bibitem{ItoTer073}
\bysame, \emph{Tridiagonal pairs and the quantum affine algebra
  ${U}_q(\widehat{\mathfrak{sl}}_2)$}, Ramanujan J. \textbf{13 (1--3)} (2007),
  39--62.

\bibitem{ItoTerinpress4}
\bysame, \emph{Finite dimensional irreducible modules for the three-point
  $\mathfrak{sl}\sb 2$ loop algebra}, Comm. Algebra \textbf{36} (2008), no.~12,
  4557--4598.

\bibitem{ItoTerinpress1}
\bysame, \emph{Distance regular graphs and the $q$-tetrahedron algebra},
  European J. Combin. \textbf{30} (2009), no.~3, 682--697.

\bibitem{ItoTerinpress2}
\bysame, \emph{Distance regular graphs of $q$-{R}acah type and the
  $q$-tetrahedron algebra}, Michigan Math. J. \textbf{58} (2009), no.~1,
  241--254.

\bibitem{ItoTer09}
\bysame, \emph{Tridiagonal pairs of $q$-{R}acah type}, J. Algebra \textbf{322}
  (2009), no.~1, 68--93.

\bibitem{ItoTerWang06}
T.~Ito, P.~Terwilliger, and C.W. Weng, \emph{The quantum algebra
  ${U}_q(\mathfrak{sl}_2)$ and its equitable presentation}, J. Algebra
  \textbf{298} (2006), 284--301.

\bibitem{Jantzen96}
J.C. Jantzen, \emph{Lectures on quantum groups}, American Mathematical Society,
  Providence, RI, 1996.

\bibitem{Miki10}
K.~Miki, \emph{Finite dimensional modules for the $q$-tetrahedron algebra},
  Osaka J. Math. \textbf{47} (2010), 559--589.

\bibitem{Ter011}
P.~Terwilliger, \emph{Two linear transformations each tridiagonal with respect
  to an eigenbasis of the other}, Linear Algebra Appl. \textbf{330} (2001),
  149--203.

\bibitem{Ter042}
\bysame, \emph{An algebraic approach to the {A}skey scheme of orthogonal
  polynomials}, {O}rthogonal polynomials and special functions, Lecture Notes
  in Math. (Berlin), vol. 1883, Springer, 2006, pp.~255--330.

\bibitem{Terinpress}
\bysame, \emph{The equitable presentation for the quantum group
  ${U}_q({\mathfrak{g}})$ associated with a symmetrizable {K}ac-{M}oody algebra
  $\mathfrak{g}$}, J. Algebra \textbf{298} (2006), no.~1, 302--319.

\bibitem{Ter13}
\bysame, \emph{Finite-dimensional irreducible ${U}_q(\mathfrak{sl}_2)$-modules
  from the equitable point of view}, Linear Algebra Appl. \textbf{439} (2013),
  no.~2, 358--400.

\end{thebibliography}

\noindent Darren Funk-Neubauer \hfill email: {\tt darren.funkneubauer@csupueblo.edu} \\                                                                                                                 \noindent Phone: (719) 549 - 2693 \hfill Fax:  (719) 549 - 2962 \\
\noindent Department of Mathematics and Physics \hfill Colorado State University - Pueblo \\
\noindent 2200 Bonforte Boulevard \hfill Pueblo, CO 81001 USA 

\end{document}